\documentclass{article}
\usepackage{amsfonts}
\usepackage{mathrsfs}

\usepackage{amssymb, latexsym}

\usepackage{graphicx}

\usepackage{amsmath}

\def\en t{{{\rm Z}\mkern-5.5mu{\rm Z}}}

\newtheorem{theorem}{Theorem}[section]

\newtheorem{conjecture}[theorem]{Conjecture}

\newtheorem{lemma}[theorem]{Lemma}

\newtheorem{remark}[theorem]{Remark}

\begin{document}

\title{\Large\bf Generalized Dyson Brownian motion, McKean-Vlasov equation and eigenvalues of random matrices}

\author{Songzi Li, \ \ Xiang-Dong Li\thanks{Research supported by NSFC No. 10971032, Key Laboratory RCSDS, CAS, No. 2008DP173182, and a
Hundred Talents Project of AMSS, CAS.}, \ \ Yong-Xiao Xie}

\maketitle

\begin{center}
\begin{minipage}{120mm}
\begin{center}{\bf Abstract}\end{center}  Using It\^o's calculus and the mass optimal transportation theory, we study the generalized Dyson Brownian motion (GDBM) and the associated McKean-Vlasov
evolution equation with an external potential $V$. Under suitable
condition on $V$, we prove the existence and uniqueness of strong
solution to SDE for GDBM. Standard argument shows that the family of
the process of empirical measures $L_N(t)$ of GDBM is tight and
every accumulative point of $L_N(t)$ in the weak convergence
topology is a weak solution of the associated  McKean-Vlasov
evolution equation, which can be realized as the gradient flow of
the Voiculescu free entropy on the Wasserstein space over
$\mathbb{R}$. Under the condition $V''\geq -K$ for some constant
$K\geq 0$, we prove that the McKean-Vlasov equation has a unique
solution $\mu(t)$ and $L_N(t)$ converges weakly to $\mu(t)$ as
$N\rightarrow \infty$. For $C^2$ convex potentials, we prove that
$\mu(t)$ converges to the equilibrium measure $\mu_V$ with respect
to the $W_2$-Wasserstein distance on $\mathscr{P}_2(\mathbb{R})$ as
$t\rightarrow \infty$. Under the uniform convexity or a modified
uniform convexity condition on $V$, we prove that $\mu(t)$ converges
to $\mu_V$ with respect to the $W_2$-Wasserstein distance on
$\mathscr{P}_2(\mathbb{R})$ with an exponential rate as
$t\rightarrow \infty$. Finally, we discuss the double-well potentials and raise some
conjectures.
\end{minipage}
\end{center}

\noindent{\bf Key words and phrases:}\ \  Generalized Dyson Brownian motion,
McKean-Vlasov equation, Johansson's theorem, gradient flow, Voiculescu free entropy, Wasserstein distance.

\section{Introduction}

\subsection{Background}

Let $V: \mathbb{R}\rightarrow \mathbb{R}^+$ be a real polynomial of even degree with
positive leading coefficient, or more general a real analytic function. Consider the following probability
measure on $\mathcal {H}_N$ (the set of
all~$N\times N$ Hermitian matrices)
\begin{eqnarray*}
d\mu_{N}(M)=\frac{1}{Z_{N}}\exp(-N{\rm Tr}V(M))dM,
\label{P1.2}
\end{eqnarray*}
where $Z_{N}$ is a normalization constant,
and if we denote $x_1, \ldots, x_N$ the eigenvalues of $M$,
$$
{\rm Tr}V(M)=\sum\limits_{i=1}^N V(x_i).$$ By \cite{BPS, Joh98}, the distribution of the eigenvalues of $M$ has the following probability density
$$\rho_{N}(x)=\frac{1}{Z_{N}}\Pi_{i<
j}|x_i-x_j|^2\exp\left(-N\sum\limits_{i=1}^NV(x_i)\right), \ \ \ x\in \mathbb{R}^N.$$

The probability distribution
$\rho_{N}(x)dx$ has a statistical mechanical interpretation: it is
the canonical Gibbs measure, at inverse temperature $\beta=2,$ of a
system of $N$ unit charges interacting through the logarithmic Coulomb
potential and confined by an external potential $NV$. From the
statistical mechanical point of view, it is natural to consider the logarithmic
Coulomb gas at arbitrary values of the inverse temperature $\beta>0$. More generally, let $V: \mathbb{R}\rightarrow \mathbb{R}$ be a polynomial of even degree or an analytic function with the growth condition
\begin{eqnarray}
V(x)\geq (1+\delta)\log(x^2+1),\ \ \ \ x\in \mathbb{R},\label{grow}
\end{eqnarray}
the $\beta$-invariant random matrix ensemble, or the so-called log-gas model,  is defined as an interacting particle system with the following probability distribution
\begin{eqnarray*}
P^N_\beta(dx_1,\ldots,dx_N)=\frac{1}{Z^\beta_N} \Pi_{i\neq
j}|x_i-x_j|^{\frac{\beta}{2}}\exp\left(-\frac{\beta N
}{2}\sum\limits_{i=1}^NV(x_i)\right)\prod_{i=1}^N dx_i,
\end{eqnarray*}
where $\beta>0$ is a parameter, called the inverse of the temperature in statistical physics.
When $\beta=1, 2, 4$, one can realize the above distribution as the distribution of the eigenvalues of $N\times N$ random matrices. More precisely, for $\beta=1, 2, 4$, if $M$ is a $N\times N$ matrix with distribution
\begin{eqnarray*}
dP^N_\beta(M)=\frac{1}{Z_N}\exp\left(-\frac{\beta N}{2}{\rm Tr} V(M)\right)dM,
\end{eqnarray*}
where $dM$ denotes the standard measure on $\mathcal{H}_N$, i.e.,
$$
dM=\prod_{1\leq i\leq N} dM_{ii}\prod_{1\leq i<j\leq N} d{\rm Re}(M_{ij})d{\rm Im}M_{ij}.$$
In particular, for $V(x)={x^2\over 2}$, we get Gaussian Symmetric matrices ensemble (GOE) for $\beta=1$, the Gaussian Hermitian matrices ensemble (GUE) for $\beta=2$, and for the Gaussian Symplectical  matrices ensemble (GSE) for $\beta=4$.

The distribution of interacting particles with general external potential $V$ and the logarithmic Coulomb interaction has  received a
lot of attention in theoretic physics in connection with the
so-called matrix models, see e.g. \cite{BZ93, FFS92}.  We would like also to mention that, Kontsevich \cite{Kon} also used the complex partition function for the matrix model with $V(x)=x^3$ on Hermitian random matrices to prove Witten's conjecture in the intersection theory of the moduli space of curves.

Suppose that $V:\mathbb{R}\rightarrow \mathbb{R}$ satisfies $(\ref{grow})$. In \cite{Joh98}, Johansson proved the following result: There is a unique equilibrium measure
$\mu_V\in\mathscr{P}(\mathbb{R})$ with compact support  such that
$$\inf\limits_{\mu\in\mathscr{P}(\mathbb{R})}\Sigma_V[\mu]=\Sigma_V[\mu_V],$$
and satisfies the Euler-Lagrange equation
\begin{eqnarray*}
{\rm P.V.} \int_{\mathbb{R}}\frac{d\mu_V(y)}{x-y}=\frac{1}{2}V'(x),\ \ \ \
x\in {\rm supp} \mu_V, \label{P1.3}
\end{eqnarray*}
where $\mathscr{P}(\mathbb{R})$ is the set of all Bore probability
measures on~$\mathbb{R}$, and $\Sigma_V$ is the following energy functional, the so-called Voiculescu free entropy functional
\begin{eqnarray*}
\Sigma_V(\mu)=-\int_\mathbb{R}\int_\mathbb{R}\log|x-y| d\mu(x)d\mu(y)+\int_\mathbb{R}V(x)d\mu(x).
\end{eqnarray*}
Moreover, as $N$ tends to infinity, the expectation of the empirical measure $L_N={1\over N}\sum\limits_{i=1}^N \delta_{\lambda_i}$  weakly converges  to $\mu_V$, i.e., $\mathbb{E}\left[L_N\right]\rightarrow \mu_V$. This recovers Wigner's famous semi-circle law \cite{W} for GUE with $V(x)={x^2\over 2}$.

Let us consider the Gaussian Hermitian Ensemble (GUE). Let $B_t$ be the $N\times N$ Hermitian matrice with
entries $B_{ij}(t)$, where $B_{ij}(t)$, $1\leq i\leq j\leq N$, are
i.i.d. complex Brownian motions with $\mathbb{E}[B_{ij}(t)]=0$ and
$\mathbb{E}[|B_{ij}|^2(t)]=t$. Let $\lambda_1(t), \ldots,
\lambda_N(t)$ be the process of the eigenvalues of $B_t$. In
\cite{Dy2}, Dyson  proved that $\lambda_i(t)$ satisfies the
following stochastic differential equations
\begin{eqnarray*} d\lambda_N^i(t)={1\over \sqrt{N}}dW^i_{t}+\frac{1}{N}\sum\limits_{j:j\neq i}
\frac{1}{\lambda^i_N(t)-\lambda^j_N(t)}dt. \label{DBM1}
\end{eqnarray*}
See also Mehta \cite{Meh}, Guionnet \cite{Gui}, Andersson-Guionnet-Zeitouni \cite{AGZ} and Tao
\cite{Tao}.

The Dyson Brownian motion $\{\lambda_i(t), i=1, \ldots, N\}$ is an
interacting $N$-particle system with the logarithmic Coulomb
interaction. It has been very useful in various branches of
mathematics and physics, including statistical physics and the
quantum chaotic systems. See the references mentioned in \cite{Meh,
Gui, AGZ, PaSh}. Instead of using Hermitian matrix-valued Brownian
motion, Chan \cite{Ch}, Rogers and Shi \cite{RS93} proved that the
eigenvalues of the $N\times N$ Hermitian Ornsetin-Uhlenbeck process
satisfies the following stochastic differential equations
\begin{eqnarray*} d\lambda_N^i(t)={1\over \sqrt{N}}dW^i_{t}+\frac{1}{N}\sum\limits_{j:j\neq i}
\frac{1}{\lambda^i_N(t)-\lambda^j_N(t)}dt-{1\over 2}\lambda_N^i(t)dt. \label{DBM2}
\end{eqnarray*}
They \cite{Ch, RS93} proved the tightness of the family of the
empirical measure $L_N(t)={1\over N}\sum\limits_{i=1}^N
\lambda_i(t)$ for the above Dyson-Ornstein-Uhlenbeck Brownian motion
$(\lambda_1(t), \ldots, \lambda_N(t))$ and proved that the limit of
any weakly convergence subsequence of $L_N(t)$ is a weak solution of
the following McKean-Vlasov equation
\begin{eqnarray*}
\frac{d}{dt}\int_{\mathbb{R}}
f(x)\mu_t(dx)=\frac{1}{2}\int\int_{\mathbb{R}^2}\frac{\partial_xf(x)-\partial_yf(y)}{x-y}\mu_t(dx)\mu_t(dy)-\frac{1}{2}\int_{\mathbb{R}}
x f'(x)\mu_t(dx), \label{DBM3}
\end{eqnarray*}
where $f$ is a test function in $C^2_b(\mathbb{R}^N)$. In
particular, the Hilbert-Stiejies transformation of $\mu(t)$, defined
by
\begin{eqnarray*}
G_t(z)=\int_{\mathbb{R}} {\mu_t(dx)\over z-x}, \label{DBM4}
\end{eqnarray*}
satisfies the following nonlinear Burgers type equation
\begin{eqnarray*}
{\partial\over \partial t}G_t(z)=\left(-G_t(z)+ z\right)
{\partial\over
\partial z}G_t(z)+ G_t(z).\label{DBM5}
\end{eqnarray*}
Showing that the above complex Burgers equation has a
unique solution, they derived that $L_N(t)$ converges weakly to
$\mu(t)$. Moreover, as $t\rightarrow \infty$, $G_t(z)$ converges to
the Hilbert-Stiejies transform of the semi-circle law, one can
therefore give a new proof of Wigner's theorem \cite{W} for the convergence
of the empirical measure of the eigenvalues of Gaussian Unitary
Ensemble (GUE) to the semi-circle law. See also Guionnet \cite{Gui}
and Andersson-Guionnet-Zeitouni \cite{AGZ} for a nice presentation
of this dynamic approach to Wigner's theorem using Dyson's Bownian
motion.

\subsection{Motivation}

The purpose of this paper is to study the generalized Dyson Brownian motion and related McKean-Vlasov equation associated with the log-gas model with non-quadratic external potential. Thus, we are working on non Gaussian type $\beta$-invariant ensembles with a general external potential $V$.
To describe our motivation, let us first introduce the generalized Dyson Brownian motions (briefly, GDBM) as follows.

Let $(W^{1},\ldots, W^{N})$ be an $N$-dimensional Brownian motion
defined on a probability space~$(\Omega,\mathbb{P})$  with a
filtration~$\mathscr{F} = \{ \mathscr{F}_t, t \geq 0\}$. Let
~$\lambda_N(0)=(\lambda^1_N(0),\ldots,\lambda^N_N(0))\in
\bigtriangleup_{N}$, where
$$\bigtriangleup_{N}=\{(x_{i})_{1\leq i\leq N}\in\mathbb{R}^{N}:
x_{1}<x_{2}<\ldots<x_{N}\}.$$ By Theorem \ref{Th1} below,  for a  wide class
of potential function $V$ with suitable condition,  the following stochastic differential
equations
\begin{eqnarray} d\lambda_N^i(t)=\sqrt{\frac{2}{\beta
N}}dW^i_{t}+\frac{1}{N}\sum\limits_{j:j\neq i}
\frac{1}{\lambda^i_N(t)-\lambda^j_N(t)}dt-{1\over 2}V'(\lambda_N^i(t))dt, \ \ \ i=1, \ldots, N, \label{SDE1}
\end{eqnarray} have a
unique strong solution ~$(\lambda_N(t))_{t\geq0}$~with initial data
~$\lambda_N(0)$~such that~$\lambda_N(t)\in\triangle_N$ for
all~$t\geq0$.

The process $(\lambda_N(t))_{t\geq0}$ defined by SDE (\ref{SDE1})
is called the generalized Dyson Brownian motion (GDBM) with potential $V$.
The GDBM is an interacting particle
system  with Hamiltonian of the form
$$H(x_1, \ldots, x_N):=-\frac{1}{2N}\sum\limits_{1\leq i\neq
j\leq N}\log |x_i-x_j|+\frac{1}{2}\sum\limits_{i=1}^NV(x_i).$$
The infinitesimal generator of these interacting particles
 is given by
\begin{eqnarray*}
\mathscr{L}^\beta_{N}f(x)=\frac{1}{\beta
N}\sum\limits_{k=1}^N\frac{\partial^2f(x)}{\partial
x_k^2}+\sum\limits_{k=1}^{N}\left({\rm P.V.}
\int_{\mathbb{R}}\frac{L_N(dy)}{x_k-y}-\frac{1}{2}V'(x_k)\right)\frac{\partial
f(x)}{\partial x_k},
\end{eqnarray*}
where  $f\in C^2(\mathbb{R}^N)$ and
$L_N=\frac{1}{N}\sum\limits_{i=1}^N\delta_{x_i}\in\mathscr{P}(\mathbb{R})$.

In particular, when $V\equiv 0$, $(\lambda_N(t))_{t\geq0}$ is the
Dyson Brownian motion introduced by F. Dyson \cite{Dy1, Dy2}, and
when $V(x)={x^2\over 2}$, it is the Dyson-Ornstein-Uhlenbeck
Brownian motion introduced by Chan \cite{Ch} and Rogers and Shi
\cite{RS93}. In fact, as the same situation as what has been done by Dyson for $V=0$, and
by Chan \cite{Ch} and Rogers-Shi \cite{RS93} for $V(x)={x^2\over
2}$,  we can prove that, in the case $\beta=1, 2, 4$, the
generalized Dyson Brownian motion can be realized as the process of
the eigenvalues of a matrix-valued diffusion process. See Remark
\ref{Th2}, Section 2.2 and Section 2.3 below.

The process of the empirical measures of GDBM  is defined by
$$L_N(t)=\frac{1}{N}\sum\limits^N_{i=1}\delta_{\lambda^i_N(t)}\in\mathscr{P}(\mathbb{R}), \ \ \ t\in [0, \infty).$$
In Section $3$, we will prove that, under suitable condition on $V$,
the family of the process of empirical measures $L_N(t)$ is tight,
and any weak convergence limit of convergent subsequence of
$L_N(t)$, denoted by $\mu(t)$, is a weak solution to the following
nonlinear McKean-Vlasov equation: for all $f\in C^2_b(\mathbb{R})$,
\begin{eqnarray*}
\frac{d}{dt}\int_{\mathbb{R}}
f(x)\mu_t(dx)=\frac{1}{2}\int\int_{\mathbb{R}^2}\frac{\partial_xf(x)-\partial_yf(y)}{x-y}\mu_t(dx)\mu_t(dy)-\frac{1}{2}\int_{\mathbb{R}}
V'(x)f'(x)\mu_t(dx).
\end{eqnarray*}
This also proves the existence of weak solutions of the
McKean-Vlasov  equation. If one can further prove the uniqueness of
the weak solution to the McKean-Vlasov equation, we can then derive
that the empirical measure $L_N(t)$ converges weakly to $\mu(t)$ as
$N\rightarrow \infty$.

\subsection{Interacting particle system}

The GDBM is a special example of interacting particle systems which
can be defined by a stochastic differential equation of the form
\begin{eqnarray*}
dx_t^i=dB_t^i-\nabla V(x_t^i)dt-{1\over N}\sum\limits_{j\neq i} \nabla W(x_t^i-x_t^j)dt,\ \ \ i=1, \ldots, N,
\end{eqnarray*}
where $V$ is a external potential, $W$ is an interacting function
between particles, and $B_t^i$ are i.i.d. Brownian motions on
$\mathbb{R}^d$. See \cite{Sph, BCCP1, CMV1, Mal, Vi1, Vi2} and
reference therein. For $N=\infty$, see e.g. \cite{KT, Os}. Let
$$\mu_t^N={1\over N}\sum\limits_{i=1}^N \delta_{x_t^i}$$
be the empirical measure of the particle system $\{x_i(t), i=1, \ldots, N\}$. Then the above SDE can be rewritten as
\begin{eqnarray*}
dx_t^i=dB_t^i-\nabla V(x_t^i)dt-\nabla (W*\mu_t^N)(x_t^i)dt.
\end{eqnarray*}
Under various assumptions which require that $V$ and $W$ are
Lipschitz, as $N\rightarrow \infty$, it has been proved that
$\mu_t^N$ converges weakly to a measure-valued process $\mu_t$ on
$\mathbb{R}$, i.e.,
\begin{eqnarray*}
\mu_t^N\rightarrow \mu_t,
\end{eqnarray*}
which is the law of a nonlinear diffusion process on $\mathbb{R}$
defined by
\begin{eqnarray*}
d\overline{X}_t=dB_t-\nabla V(\overline{X}_t)dt-\nabla
(W*\mu_t)(\overline{X}_t)dt.
\end{eqnarray*}
where $B_t$ is a Brownian motion on $\mathbb{R}^d$. See \cite{Mal}.
Suppose that $\mu_t<<dx$, then the density function $u={d\mu_t\over
dx}$ satisfies the nonlinear McKean-Vlasov equation (called also the
nonlinear Fokker-Planck equation in \cite{ BCCP1}  etc.)
\begin{eqnarray}
\partial_t u={\rm div}(u\nabla(\log u+V+W*u)).\label{NMV}
\end{eqnarray}

In \cite{BCCP1}, Benedetto, Caglioti, Carrillo and Pulvirenti
developed the $L^1$-theory of the nonlinear McKean-Vlasov  equation
$(\ref{NMV})$ in one-dimensional granular media. Assuming that the
interaction function $W$ and the potential $V$ are  Lipschitz and
convex functions, they proved that the free energy functional
\begin{eqnarray*}
F(u)=\int_{\mathbb{R}} u\log udx+\int_{\mathbb{R}} V
udx+\int\int_{\mathbb{R}^2} W(x-y)u(x)u(y)dxdy
\end{eqnarray*}
has a unique minimum $u_\infty$, and it holds that
\begin{eqnarray*}
\|u_t-u_\infty\|_{L^1}\rightarrow 0 \ \ {\rm as}\ \ t\rightarrow \infty.
\end{eqnarray*}
Moreover, under the Lipschitz condition on $V$ and $W$, they proved that there exists a constant $C>0$ such that
\begin{eqnarray*}
W_1(\mu_t, \nu_t)\leq e^{Ct}W_1(\mu_0, \nu_0),
\end{eqnarray*}
where $\mu_t$ and $\nu_t$ are solutions to the  nonlinear
McKean-Vlasov equation  $(\ref{NMV})$  with initial dates $\mu_0$
and $\nu_0$ respectively. See also \cite{Sph, Lan} and references
therein.

However, the logarithmic Coulomb interaction function  appeared in
the distribution of the eigenvalues of $\beta$-invariant random
matrices ensemble, i.e.,
$$W(x)=-\log |x|$$
is not a Lipschitz function on $\mathbb{R}$. Thus, it does not
satisfy the Lipschitz condition required in \cite{Mal, BCCP1},  and
the $L^1$-theory of Benedetto et al.\cite{BCCP1} does not apply
directly to $(\ref{NMV})$ associated to the GDBM, even in the case
of Gaussian ensembles, i.e., $V(x)={x^2\over 2}$. As far as we know, the
nonlinear McKean-Vlasov equation $(\ref{NMV})$ with the logarithmic
Coulomb interaction $W(x)=-\log|x|$ and general external potential
$V$ has not been well studied in the literature. In particular, even
though the existence of weak solution of $(\ref{NMV})$ can be easily
derived from the McKean-Vlasov limit of any convergent subsequence
of the empirical measure of the generalized Dyson Brownian motion,
the problem of the uniqueness and the longtime asymptotic behavior
of the solutions of the nonlinear McKean-Vlasov equation
$(\ref{NMV})$  remain as an open problem in the literature.
\medskip

\subsection{Otto's calculus and the gradient flow of the Voiculescu entropy}

In \cite{CMV1}, Carrillo, McCann and Villani used Otto's  infinite
dimensional differential calculus \cite{Ot} on the  Wasserstein
space over $\mathbb{R}^n$ to study the convergence rate problem of
the nonlinear McKean-Vlasov equation $(\ref{NMV})$ in the granular
media. Under some growth and smoothness assumptions of the
interaction potential $W$ and the external potential $V$, they
proved that, if $\nabla^2W\geq 0$ and $\nabla^2V\geq \lambda$, or if
$\nabla^2W\geq \lambda$ and $\nabla^2V\geq 0$, where $\lambda>0$ is
a constant, then  the solution $u(t)$ of $(\ref{NMV})$ converges to
the equilibrium $u_\infty$ in the $W_2$-Wassersten distance with an
exponential rate, i.e.,
\begin{eqnarray*}
W_2(u_t, u_\infty)=O(e^{-\lambda t}).
\end{eqnarray*}

Note that, for  $W(x)=-\log |x|$, we have
\begin{eqnarray*}
\nabla^2 W(x)={1\over |x|^2}, \ \ \ \ x\neq 0.
\end{eqnarray*}
Thus, the logarithmic Coulomb interaction potential has a strong
convexity near its singularity point $x=0$. This suggests us to
adopt the  infinite dimensional calculus on the Wasserstein space
initiated by Otto \cite{Ot} and  developed by Carrillo, McCann and
Villani \cite{CMV1} to study the uniqueness problem and the long
time asymptotic behavior of the solution of the McKean-Vlasov
evolution equation.

Let $V$ be a $C^1$ function on $\mathbb{R}$. According to Voiculescu \cite{Voi1}, we introduce the free entropy as follows
\begin{eqnarray*}
\Sigma_V(\mu)=-\int\int_{\mathbb{R}^2} \log|x-y|d\mu(x)d\mu(y)+\int_{\mathbb{R}}V(x)d\mu(x).
\end{eqnarray*}It has the following electrostatic interpretation: Suppose that electrons are distributed on the one-dimensional space $\mathbb{R}$ in the presence of an external field with potential $V$. For a probability distribution $\mu$ of electrons, the electrostatic repulsion is given by
\begin{eqnarray*}
\int\int_{\mathbb{R}^2}\log|x-y|^{-1}d\mu(x)d\mu(y),
\end{eqnarray*}
and the energy from the external field is $\int_{\mathbb{R}} V(x)d\mu(x)$. Thus, the Voicuselecu free entropy is the total energy of electrons in an external field with potential $V$. When $V$ satisfies the growth condition $(\ref{grow})$, it is shown in \cite{Joh98} that there exists a unique minimizer of the Voiculescu free entropy
\begin{eqnarray*}
\mu_V={\rm arg min}_{\mu\in \mathcal{P}(\mathbb{R})} \Sigma_V(\mu),
\end{eqnarray*}
called the equilibrium distribution of $\Sigma_V$.
The relative free entropy is defined as follows
\begin{eqnarray*}
\Sigma_V(\mu|\mu_V)=\Sigma_V(\mu)-\Sigma_V(\mu_V).
\end{eqnarray*}

Following Voiculescu \cite{Voi1} and Biane
\cite{Bian03}, the relative free Fisher information is defined as
follows
\begin{eqnarray*}
{\rm I}_V(\mu)=\int_{\mathbb{R}}\left(H\mu(x)-{1\over 2}V'(x)\right)^2d\mu(x),
\end{eqnarray*}
where
\begin{eqnarray*}
H\mu:={\rm PV}. \int_{\mathbb{R}} {d\mu(y)\over x-y}.
\end{eqnarray*}
is the Hilbert transform of $\mu$.

Note that, as the equilibrium measure $\mu_V$ satisfies the equation
\begin{eqnarray*}
H\mu_V(x)={1\over 2}V'(x),\ \ \ \ \forall x\in \mathbb{R},
\end{eqnarray*}
we have
\begin{eqnarray*}
{\rm I}_V(\mu_V)=0.
\end{eqnarray*}

Inspired by Otto \cite{Ot}, Carrillo-McCann-Villani \cite{CMV1}, Villani \cite{Vi1, Vi2} and  Biane \cite{Bian03}, we have the following result, which has been known by Biane and Speicher \cite{Bian-Sp01}.\\

\begin{theorem}\label{Th0} (Biane and Speicher \cite{Bian-Sp01}) Let $V\in C^2(\mathbb{R})$. Then the McKean-Vlasov equation $(\ref{DBM7})$ is the gradient flow of the Voiculescu free entropy $\Sigma_V$ on the Wasserstein space $\mathscr{P}_2(\mathbb{R})$ over $\mathbb{R}$ equipped with Otto's infinite dimensional Riemannian structure.
\end{theorem}

For the definition of the Otto's infinite dimensional Riemannian
structure on the Wasserstein space $\mathscr{P}_2(\mathbb{R})$, see
Section $4$.

\subsection{Main results}

We  now state the main results of this paper. Our first result gives the existence and uniqueness of the strong solution to the SDE of the generalized Dyson Brownian motion under reasonable condition on the external potential.

\begin{theorem}\label{Th1}
Let~$(W^{1},\ldots,W^{N})$~be an $N$-dimensional Brownian motion in
a probability space~$(\Omega,\mathbb{P})$~equipped with a
filtration~$\mathscr{F} = \{ \mathscr{F}_t, t \geq 0\}$.  Let $V\in C^2(\mathbb{R})$ be a function satisfying the growth condition $(\ref{grow})$ and the following conditions\\
(i) For all $R>0$, there is $K_R>0,$~such that for all
~$x,y\in\mathbb{R}$ with $|x|, |y|\leq R$, $$(x-y)(V'(x)-V'(y))\geq -K_R|x-y|^2,$$~\\
(ii) For some~$\gamma>0,$~such that
\begin{eqnarray}
-xV'(x)\leq \gamma(1+|x|^2), \ \ \forall \ x\in
\mathbb{R}.\label{Cond-LX1}
\end{eqnarray}
Then, for any $\beta\geq 1$, and for any given
$\lambda_N(0)=(\lambda^1_N(0),\ldots,\lambda^N_N(0))\in
\bigtriangleup_{N}$, there exists a unique strong solution
$(\lambda_N(t))_{t\geq0}$ with infinite lifetime to SDE
$(\ref{SDE1})$, with initial value $\lambda_N(0)$ and such that
$\lambda_N(t)\in\triangle_N$~for all~$t\geq0$.
\end{theorem}

The second result of this paper is the following result concerning the existence and uniqueness of the weak solution to the McKean-Vlasov equation, as well as the convergence of the empirical measure $L_N(t)$ towards the solution of the McKean-Vlasov equation.

\begin{theorem}\label{Th4} (i)Suppose that $V$ be a $C^2$ function satisfying the same condition as in Theorem \ref{Th1}. Suppose that
\begin{eqnarray*}
\sup\limits_{N\geq 0}\int_{\mathbb{R}}\log(x^2+1)dL_N(0)<\infty,
\end{eqnarray*}
and
$$L_N(0)=\frac{1}{N}\sum\limits_{k=1}^N\delta_{\lambda^k_N(0)}\rightarrow\mu\in\mathscr{P}(\mathbb{R}) \ \ \ {\rm as}\ \ N\rightarrow\infty.$$ Then, the family $\{L_N(t), t\in [0, T]]\}$ is tight, and the limit of any weakly convergent subsequence of $\{L_N(t), t\in [0, T]]\}$ is a weak solution of the McKean-Vlasov equation, i.e., for all $f\in C_b^2(\mathbb{R})$, $t\in [0, T]$,
\begin{eqnarray}
\frac{d}{dt}\int_{\mathbb{R}}
f(x)\mu_t(dx)=\frac{1}{2}\int\int_{\mathbb{R}^2}\frac{\partial_xf(x)-\partial_yf(y)}{x-y}\mu_t(dx)\mu_t(dy)-\frac{1}{2}\int_{\mathbb{R}}
V'(x)f'(x)\mu_t(dx). \label{DBM7}
\end{eqnarray}
Equivalently, the probability density of $\mu_t$ with respect to the Lebesgue measure on $\mathbb{R}$ satisfies the nonlinear Fokker-Planck equation
\begin{eqnarray}
{\partial \rho_t\over \partial t}={\partial \over \partial
x}\left(\rho_t\left({1\over 2}V'-H\rho_t\right)\right),\label{NFK1}
\end{eqnarray}
where
$$
{\rm H}\rho(x)={\rm P.V.}\int_{\mathbb{R}}{\rho(y)\over x-y}dy$$ is
the Hilbert transform of $\rho$.\\
(ii) Suppose that there exists a constant $K\in \mathbb{R}$ such that
$$V''(x)\geq K, \ \ \ \ x\in \mathbb{R}.$$
Let $\mu_i(t)$ be two solutions of the McKean-Vlasov equation
$(\ref{DBM7})$ with initial datas $\mu_i(0)$, $i=1, 2$. Then for all
$t>0$, we have
\begin{eqnarray*}
W_2(\mu_1(t), \mu_2(t))\leq e^{-Kt}W_2(\mu_1(0), \mu_2(0)).
\end{eqnarray*}
In particular, the McKean-Vlasov equation $(\ref{DBM7})$ has a
unique solution.\\
 (iii) Let $V$ be a $C^2$ function satisfying the same condition as in Theorem \ref{Th1} and $V''\geq K$ for some constant $K\in \mathbb{R}$. Then the empirical measure $L_N(t)$ weakly converges to the unique solution $\mu(t)$ of the McKean-Vlasov equation
$(\ref{DBM7})$.
\end{theorem}

The following result gives us the longtime asymptotic behavior and the convergence rate of the weak solution of the McKean-Vlasov equation to the equilibrium measure of the Voiculescu free entropy for $C^2$-convex potentials.

\begin{theorem} \label{Th5}  (i) Let $V$ be a $C^2$-convex potential. Then $\mu(t)$ converges to $\mu_V$ with respect to the Wasserstein distance in $\mathscr{P}_2(\mathbb{R})$, i.e.,
\begin{eqnarray*}
W_2(\mu(t), \mu_V)\rightarrow 0 \ \ \ {\rm as}\ \  t\rightarrow \infty.
\end{eqnarray*}
(ii) Suppose that there exists a constant $K\in \mathbb{R}$ such that
$$V''(x)\geq K, \ \ \ \forall x\in \mathbb{R}.$$
Then for all $t>0$, we have
\begin{eqnarray*}
\Sigma_V(\mu_t|\mu_V)&\leq& e^{-2Kt}
\Sigma_V(\mu_0|\mu_V),\\
W_2(\mu_t, \mu_V)&\leq& e^{-Kt}W_2(\mu_0, \mu_V).
\end{eqnarray*}
In particular, if $K>0$, then $\mu_t$ converges to $\mu_V$ with respect to the
$W_2$-Wasserstein distance with the exponential rate $K$.\\
(iii) Suppose that $V$ is a $C^2$, convex and there exists a constant $r>0$ such that
$$V''(x)\geq K>0, \ \ \ \ |x|\geq r.$$
Then there exist constants $C_1>0$ and $C_2>0$ such that $\mu(t)$
converges to $\mu_V$ with respect to the $W_2$-Wasserstein distance  with an
exponential rate
\begin{eqnarray*}
W^2_2(\mu_t, \mu_V)\leq {e^{-C_1t}\over C_2}\Sigma_V(\mu_0|\mu), \ \ \ \ \ t>0.
\end{eqnarray*}
\end{theorem}

\subsection{Remarks}

\begin{remark}\label{rem1} $(i)$ In \cite{RS93},  Rogers and Shi proved (even though did not state) the existence and uniqueness of
the generalized Dyson Brownian motion under the condition that the
potential $V$ satisfies
\begin{eqnarray}
-xV'(x)\leq \gamma, \ \ \forall \ x\in \mathbb{R}.\label{Cond-RS}
\end{eqnarray}
It is clear that our condition $(\ref{Cond-LX1})$ on $V$ is weaker
than Rogers-Shi's condition $(\ref{Cond-RS})$. To prove Theorem \ref{Th1}, we need to use the non-explosion criterion of the Bessell type SDE on the half line $\mathbb{R}^+$ (cf. \cite{IW} p. 235-237), which was kindly pointed to us by Yves Le Jan in June 2011. \\
$(ii)$ The conclusion of Theorem \ref{Th2} says that, for any given
$\lambda_N^1(0)<\lambda_N^2(0)<\ldots<\lambda_N^{N-1}(0)<\lambda_N^N(0)$,
SDE $(\ref{SDE1})$ admits a unique strong solution with infinite
lifetime such that
$\lambda_N^1(t)<\lambda_N^2(t)<\ldots<\lambda_N^{N-1}(t)<\lambda_N^N(t)$.
Therefore, the generalized Dyson Brownian particle system
$(\lambda_N^1(t), \ldots, \lambda_N^N(t))$ does not self-intersect
for all time $t>0$.
\end{remark}

\begin{remark}\label{Th2}
Let~$\beta=1, 2, 4$.
Let $V$ be an analytic function satisfying the assumptions of Theorem
\ref{Th1}. Following an original idea due to Dyson \cite{Dy2} and developed by other authors \cite{Ch, RS93, AGZ, Gui} for the special case $V(x)=\theta x^2$, we can prove that the generalized Dyson Brownian motion
$\lambda_N^1(t)\leq\ldots  \leq\lambda_N^N(t)$ can be realized as
the eigenvalues of a matrix-valued diffusion process $X_t$, which
satisfies the following SDE
\begin{eqnarray}
dX_t=\sqrt{2\over \beta N}dW_t-{1\over 2}\nabla{\rm Tr}V(X_t)dt,
\end{eqnarray}
where $B_t$ is the standard Brownian motion on
$\mathcal{H}_N^\beta$ (see Section $2.1$), and  $\nabla$ denotes the gradient operator on
$\mathcal{H}_N^\beta$. Moreover, $\lambda_N(t)=(\lambda_1(t), \ldots, \lambda_N(t)) $ is a
$\Delta_N$-valued semi-martingale. See Theorem \ref{LLX1} and Theorem \ref{Th-LX2}.
\end{remark}

\begin{remark}\label{rem2} In \cite{Bian03},
Biane pointed out that there exists non-convex
potentials such that the equation $H\rho_i(x)={1\over 2}V'(x)$ on the supports of
$\mu_i$ holds for distinct measures $\mu_i=\rho_i(x)dx$. The simplest
example due to Biane \cite{Bian03} is a two-well potential $V$
satisfying $V(x)={1\over 2}(x-a_i)^2$ in $[a_i-2, a_i+2]$ for points
$a_1, a_2$ with $|a_1-a_2|>4$, then the semi-circular measures
centered on $a_1$ and $a_2$ respectively satisfy $Hp_i(x)=V'(x)$ on
their supports respectively. By Theorem \ref{Th1} (due to Biane and
Speicher \cite{Bian-Sp01}), the McKean-Vlasov limit
$\mu(t)=\rho_t(x)dx$ of $L_N(t)$ satisfies the free Fokker-Planck
equation
\begin{eqnarray*}
{\partial \rho_t\over \partial t}={\partial \over \partial
x}\left(\rho_t\left({1\over 2}V'-H\rho_t\right)\right).
\end{eqnarray*}
It is natural to ask the question whether $\mu(t)$ converges
to $\mu_V$ in the weak convergence topology or with respect to the $W_2$-Wasserstein distance (or the $W_p$-Wasserstein distance for some $p\geq 1$)
for non-convex potentials $V$. In
\cite{Bian-Sp01}, it was pointed out that this cannot be true in
general. Indeed, as $\mu_t$ satisfies the gradient flow of the
Voiculescu free entropy $\Sigma_V$ on $\mathscr{P}(\mathbb{R})$,
$\mu_t$ may converge to a local minimizer of $\Sigma_V$ which is not
necessary the global minimizer $\mu_V$. See also Section $6$ and in particular Conjecture \ref{conj}.
\end{remark}

The rest of this paper is organized as follows. In Section $2$, we
prove Theorem \ref{Th1}. In Section $3$, we
prove Theorem \ref{Th4} (i). In Section $4$, we prove Theorem \ref{Th5}. In Section $5$,  we prove Theorem \ref{Th4} (ii) and (iii).  Finally, we study the double-well potentials and raise some conjectures in
Section $6$.

\medskip

\noindent{\bf Acknowledgement}. \  During the preparation of this
paper, the authors have benefited from some discussions with G. Akemann, F. G\"otze, K. Kuwae, Y. Le Jan,  Y. Liu, T. Lyons, S.-G. Peng,
M. R\"ockner, W. Sun and M. Venker. We would
like to express our gratitude to all of them for their interests and
for useful comments.

\section{Generalized Dyson's Brownian motion (GDBM)}

\subsection{Proof of Theorem \ref{Th1}}

For $T\in\mathbb{R}^{+},$ we denote by $\mathcal
{C}([0,T],\mathbb{R}^N)$ the space of continuous process from
$[0,T]$ to $\mathbb{R}^N$ and $\mathscr{P}(\mathcal
{C}([0,T],\triangle_N))$ is the set of all probability measures on
$\mathcal {C}([0,T],\triangle_N).$

\textbf{(1)} Fix $R>0$. Define the truncated Dyson Brownian motion
by
\begin{equation}
 d\lambda^i_{N,R}(t)=\sqrt{\frac{2}{\beta
N}}dW^i_{t}+\frac{1}{N}\sum\limits_{j:j\neq i}
\phi_R(\lambda^i_N(t)-\lambda^j_N(t))dt-\frac{1}{2}V'(\lambda^i_N(t))dt,\label{E2.1}
\end{equation}
with $\lambda^i_{N,R}(0)=\lambda^i_N(0)$~for~$1\leq i\leq N$, where
~$\phi_R(x)=
        x^{-1}$~  if~$|x|\geq R^{-1}$,~
        and~$\phi_R(x)=R^2x$ if $|x|< R^{-1}$. By Theorem 3.1.1 in \cite{Ro}, since ~$\phi_R$~is uniformly Lipschitz and $V$ satisfies the
conditions $(i)$ and $(ii)$, SDE (\ref{E2.1}) has a unique strong
solution, adapted to the filtration $\mathscr{F}$. Let
\begin{eqnarray*}
\tau_R:=\inf\{t:\min_{i\neq
j}\mid\lambda^i_{N,R}(t)-\lambda^j_{N,R}(t)\mid<R^{-1}\}.
\end{eqnarray*}
Then $\tau_R$ is monotone increasing in $R$ and
\begin{eqnarray}
\lambda_{N,R}(t)=\lambda_{N,R'}(t) \ \ \ {\rm for\ all} \ \ t\leq
\tau_R\ \ {\rm
 and} \ \  R<R'.
\end{eqnarray}

 \textbf{(2)} We now construct a solution to SDE
$(\ref{SDE1})$ by taking ~$\lambda_N(t)=\lambda_{N,R}(t)$~on the
event $[\tau_R>t]=\{|\lambda^i_N(s)-\lambda^j_N(s)|\geq R^{-1},
~\forall s\leq t, \forall i\neq j\}$. To obtain a global solution
$\{\lambda_N(t), t\in \mathbb{R}^+\}$ to SDE $(\ref{SDE1})$, we need
only to prove that $\tau_R$ tends to infinity almost surely when
$R\rightarrow \infty$.

To this end, let us consider the Lyapounov function~
\begin{eqnarray}
f(x_1,\ldots,x_N)=\frac{1}{N}\sum\limits_{i=1}^N
V(x_i)-\frac{1}{N^2}\sum\limits_{i\neq j}\log|x_i-x_j|.\label{C0}
\end{eqnarray}
Using $\log|x-y|\leq \log(|x|+1)+\log(|y|+1)$ and $(\ref{grow})$, we
can derive that
\begin{eqnarray}
f(x_1, \ldots, x_N)\geq -(1+\delta)\log 2.\label{C1}
\end{eqnarray}
and $V(x)-2\log(|x|+1)\geq C$~for  $C:=-2(1+\delta)\log 2$.
Moreover, for all~$i\neq j$, it holds
\begin{eqnarray}
-\frac{1}{N^2}\log|x_i-x_j|\leq f(x_1,\ldots,x_N)-C.\label{C2}
\end{eqnarray}

For any $M>0$, set ~$$T_M=\inf\{t\geq0:f(\lambda_N(t))\geq M\}.$$
Then, ~$T_M$~is a stopping time. Moreover, on~$\{T_M\geq T\}$, we
have
$$|\lambda_N^i(t)-\lambda_N^j(t)|\geq R^{-1}:=e^{N^2(-M+C)}, \ \
\forall\ t\leq T.$$ Thus, on the event $[T\leq T_M]$, $\{\lambda_{N,
R}(t), t\leq T\}$ is an adapted solution to SDE $(\ref{SDE1})$. It
remains to prove that for all $t\geq 0$, we have
\begin{eqnarray}
\mathbb{P}(\exists M\in \mathbb{N}: T_M\geq t)=1. \label{C3}
\end{eqnarray}

 \textbf{(3)}
To prove $(\ref{C3})$, we need only to prove that, for all $K>0$, we
have
\begin{eqnarray}
\mathbb{P}(\exists M\in \mathbb{N}: T_M\wedge \zeta_K \geq t)=1,
\label{C4}
\end{eqnarray}
where
$$\zeta_K=\inf\limits\{t\geq 0: \lambda_N^i(t)\notin [-K,
K],\ {\rm for\ some}\ i=1, \ldots, N\}.$$ By $(\ref{C0})$ and
$(\ref{C2})$, to prove $(\ref{C4})$, it is equivalent to show that
almost surely  $\lambda_N^i(t)$ and $\lambda_N^j(t)$ never collide
up to $\zeta_M$. We shall prove this claim below.

\textbf{(4)}
 By It\^{o}'s formula, for all $f\in C^2(\mathbb{R})$,
we have
\begin{eqnarray*}
df(\lambda_N(t))&=&\frac{1}{N^2}\sum\limits_{i=1}^N\left(V'(\lambda_N^i(t))-\frac{1}{N}\sum_{k\neq
i}\frac{1}{\lambda_N^i(t)-\lambda_N^k(t)}\right)\left(\sum\limits_{j\neq
i}\frac{1}{\lambda_N^i(t)-\lambda_N^j(t)}\right)dt\nonumber\\
&&+\frac{1}{2N}\sum\limits_{i=1}^N\left(-|V'(\lambda_N^i(t))|^2+\frac{1}{N}\sum\limits_{j\neq
i}\frac{V'(\lambda_N^i(t))}{\lambda_N^i(t)-\lambda_N^j(t)}\right)dt\nonumber\\
&&+\frac{1}{\beta
N^2}\sum\limits_{i=1}^N\left(V''(\lambda_N^i(t))+\frac{1}{N}\sum\limits_{j\neq
i}\frac{1}{(\lambda_N^i(t)-\lambda_N^j(t))^2}\right)dt+dM_N(t),\label{E2.2}
\end{eqnarray*}
where $M_N$ is the following local martingale
$$dM_N(t)=\frac{2^{\frac{1}{2}}}{\beta^{\frac{1}{2}}
N^{\frac{3}{2}}}\sum\limits_{i=1}^N \left(
V'(\lambda_N^i(t))-\frac{1}{N}\sum\limits_{k\neq
i}\frac{1}{\lambda_N^i(t)-\lambda_N^k(t)} \right)dW_t^i.$$ By
\cite{Gui}, we have
\begin{eqnarray*}
\sum\limits_{i=1}^N\left[\sum_{k\neq
i}\frac{1}{\lambda_N^i(t)-\lambda_N^k(t))}\sum_{j\neq
i}\frac{1}{\lambda_N^i(t)-\lambda_N^j(t)} -\sum_{j\neq
i}\frac{1}{(\lambda_N^i(t)-\lambda_N^j(t))^2}\right]=0.
\end{eqnarray*}
Thus
\begin{eqnarray*}
df(\lambda_N(t))&=&dM_N(t)+\frac{1}{N^3}\left(\frac{1}{\beta}-1\right)\sum\limits_{k\neq
i}\frac{1}{(\lambda_N^i(t)-\lambda_N^k(t))^2}dt-\frac{1}{2N}\sum\limits_{i=1}^N|V'(\lambda_N^i(t))|^2dt\\
& & +\frac{1}{N^2}\left({1\over \beta}\sum\limits_{i=1}^N
V''(\lambda_N^i(t))+\frac{3}{2}\sum\limits_{j\neq
i}\frac{V'(\lambda_N^i(t))-V'(\lambda_N^j(t))}{\lambda_N^i(t)-\lambda_N^j(t)}\right)dt\\
&=&dM_N(t)+A_N(t)dt,
\end{eqnarray*}
where, as $\beta\geq 1$, it holds
\begin{eqnarray*}
A_N(t)\leq \frac{1}{N^2}\left({1\over \beta}\sum\limits_{i=1}^N
V''(\lambda_N^i(t))+\frac{3}{2}\sum\limits_{j\neq
i}\frac{V'(\lambda_N^i(t))-V'(\lambda_N^j(t))}{\lambda_N^i(t)-\lambda_N^j(t)}\right)dt.
\end{eqnarray*}

We now prove the following lemma, which is a stronger version of
the claim in $\textbf{(3)}$.

\begin{lemma}Suppose that the processes~$(\lambda_N(t))_{t\geq 0}$~ defined by
$(\ref{SDE1})$, at least up to the stopping time
~$$\zeta=\inf\{t:\lambda_N^i(t)=\lambda_N^j(t)~\exists~i\neq j~{\rm
or}~\lambda_N^i(t)=\infty~\exists~i\}.$$ Then
$$\mathbb{P}(\zeta=\infty)=1.$$
\end{lemma}
{\it Proof}. Fix $K>0$, $T>0$ and $R>0$ such that $\lambda_N^i(0)\in
[-K, K]$ for all $i=1, \ldots, N$, and
$|\lambda_N^i(0)-\lambda_N^j(0)|\geq R^{-1}$ for all $i\neq j$, $i,
j=1, \ldots, N$.  Let $C_1(K)\geq 0$ be such that $V''(x)\leq
C_1(K)$ for all $x\in [-K, K]$. Then, $A_N(t)\leq C_1(K)$ and
$\{f(\lambda_N(t\wedge \zeta_K)-C_1(K)(t\wedge \zeta_K),\ \ t\in [0,
T]\}$ is a super-martingale.

Let $A:=\{\tau_R\leq \zeta_K,\tau_R\leq T\}$,
and $C_2(K):= \inf\limits\{V(x):|x|\leq K\}$. ~ Then
\begin{eqnarray*}
&&f(\lambda_N(0))+TC_1(K)\geq\mathbb{E}(f(\lambda_N(T\wedge
\zeta_K\wedge \tau_R)))\\
&=&\mathbb{E}\left(f(\lambda_N(\tau_R))1_A\right)+\mathbb{E}\left(f(\lambda_N(T\wedge
\zeta_K))1_{A^c}\right)\\
&\geq&\left(\frac{1}{N^2}\log R-\frac{1}{N^2}(N^2-N-1)\log(2K)+C_2(K)\right)\mathbb{P}(A)\\
& &\ \ \ \ \ \ \ \ \ \ +\left(-\frac{1}{N^2}(N^2-N)\log(2K)+C_2(K)\right)\mathbb{P}(A^c)\\
&=&\left(\frac{1}{N^2}\log
R+\frac{1}{N^2}\log(2K)\right)\mathbb{P}(A)-\frac{1}{N^2}(N^2-N)\log(2K)+C_2(K),
\end{eqnarray*}
whence
\begin{eqnarray*}
\mathbb{P}(\tau_R\leq \zeta_K,\tau_R\leq T)\leq
\frac{N^2(f(\lambda_N(0))+TC_1(K))+N(N-1)\log(2K)-C_2(K)}{\log(2K)+\log
R}.
\end{eqnarray*}
Taking $R\rightarrow \infty$,  for all $K$ and $T$, we have
$\mathbb{P}(\tau_\infty\leq \zeta_K\wedge T)=0.$ Letting $T$ tend to
infinity, we obtain $\mathbb{P}(\tau_\infty<\zeta_K)=0,$ which yields
that $\sum\limits_{i\neq
j}\log|\lambda^i_N(t)-\lambda^j_N(t)|>-\infty$ almost surely for all
$t<\zeta_K$. Moreover, letting $K$ tend to infinity, we get
$$\mathbb{P}(\tau_\infty<\zeta)=0,$$
where $\zeta:=\lim\limits_{K\rightarrow \infty}\zeta_K=\inf\limits\{t>0: f(\lambda_N(t))=\infty\}$ is the
explosion time of $f(\lambda_N(t))$. This means that
the particles $\lambda_N^1(t), \ldots, \lambda_N^N(t)$, could not
collide before the explosion.

We now prove that~$\zeta$~is infinity almost surely. To this end,
let~$$R_t=\frac{1}{2N}\sum\limits_{j=1}^N\lambda_N^j(t)^2=\langle
L_N(t),f\rangle,$$ where $f(x)={x^2\over 2}$.~ Now
$\sum\limits_{1\leq j\neq r\leq
N}\frac{\lambda_N^j(t)}{\lambda_N^j(t)-\lambda_N^r(t)}={N(N-1)\over
2}$. By It\^{o}'s formula, we have
\begin{eqnarray*}
dR_t&=&\sqrt{\frac{2}{\beta
N}}\frac{1}{N}\sum\limits_{j=1}^N\lambda_N^j(t)dW_t^j+\frac{1}{\beta
N}dt+\frac{1}{N^2}\sum\limits_{j\neq
r}\frac{\lambda_N^j(t)}{\lambda_N^j(t)-\lambda_N^r(t)}dt-\frac{1}{2}\langle
L_N(t),xV'(x)\rangle dt\\
&=&\sqrt{\frac{2}{\beta
N}}\frac{1}{N}\sum\limits_{j=1}^N\lambda_N^j(t)dW_t^j+\left(\frac{1}{\beta
N}+\frac{N-1}{2N}-\frac{1}{2}\langle L_N(t),xV'(x)\rangle\right) dt.
\end{eqnarray*}
By introducing a new Brownian motion~$B,$~we have
~$$dR_t=\sqrt{\frac{2}{\beta
N}}\frac{1}{N}\sqrt{2NR_t}dB_t+\left(\frac{1}{\beta
N}+\frac{N-1}{2N}-\frac{1}{2}\langle L_N(t),
xV'(x)\rangle\right)dt.$$ Under the assumption $(\ref{Cond-LX1})$,
and using the comparison theorem of one dimensional SDEs, cf.
\cite{IW}, we can derive that
$$R_t\leq
R'_t,\ \ \ \forall \ t\geq 0,~~{\rm a.s},$$ where~$R'$~solves
~$$dR'_t={2\over N}\sqrt{R'_t\over \beta}dB_t+\left(\frac{1}{\beta
N}+\frac{N-1}{2N}+\frac{1}{2}\gamma+\gamma R'_t\right)dt,$$ with
$R'_0=0.$ Moreover, by Example 8.2, $(8.12)$ in Ikeda and Watanabe
\cite{IW} (p. 235-237), the process $R'$ never explodes. So the
process $R$~cannot explode in finite time.

\textbf{(5)} By the continuity of the trjectory of $\lambda_N(t)$,
it is easy to see that $\lambda_N(t)\in\Delta_N$ for all $t\geq 0$.
The proof of Theorem \ref{Th1} is completed. \hfill $\square$

\medskip
\subsection{From matrix diffusions to GDBM}

The following result indicates a way to introduce generalized Dyson Brownian motion as the eigenvalues process of a matrix valued diffusion process.
To simplify notation, let $\mathcal{H}_N^{\beta}$ be the ensemble of $N\times N$ matrices: for $\beta=1$, it denotes the $N\times N$ real symmetric ensemble, for $\beta=2$, it denotes the $N\times N$ Hermitian complex ensemble, and for $\beta=4$, it denotes the $N\times N$ symplectic ensemble.

\begin{theorem}\label{LLX1} Let $\beta=1, 2, 4$, and $V:\mathbb{R}\rightarrow \mathbb{R}$
be a real analytic function. Let $X_t$ be a $\mathcal{H}_N^{\beta}$-valued
diffusion process defined by
\begin{eqnarray}
dX_t=\sqrt{2\over \beta N}dB_t-{1\over 2}\nabla{\rm
Tr}V(X_t)dt,\label{DBM1}
\end{eqnarray}
where $B_t$ is the standard Brownian motion on
$\mathcal{H}_N^\beta$, $\nabla$ denotes the gradient operator on
$\mathcal{H}_N^\beta$. Let $\lambda_N(t)=(\lambda_N^1(t), \ldots,
\lambda_N^N(t))$ be the eigenvalues of $X_N(t)$. Then,
$\lambda_N(t)$ satisfies the SDE $(\ref{SDE1})$ for the generalized
Dyson Brownian motion with $\beta=1, 2, 4$.
\end{theorem}
{\it Proof}. We only prove Theorem \ref{LLX1} for $\beta=2$. Note that, for analytic function
$V(x)=\sum\limits_{k=0}^\infty a_kx^k$ on $\mathbb{R}$, we have
\begin{eqnarray*}
\nabla {\rm Tr}V(X)=\sum\limits_{k=1}^\infty
ka_kX^{k-1}=V'(X).
\end{eqnarray*}
Hence the SDE $(\ref{DBM1})$ for $X_t$ can be written as follows
\begin{eqnarray} dX_t={1\over \sqrt{N}}dB_t-{1\over 2}V'(X)dt.\label{DBM2}
\end{eqnarray}

Let $x_1(t)\leq \ldots \leq x_N(t)$ be the ordered eigenvalues
of $X_t$. Let
$$f: X\rightarrow D={\rm diag}(x_1, \ldots, x_N)$$
be the matrix transformation such that $X=UDU^*$, where $U=(u_1,
\ldots, u_N)$ is an unitary matrix. Equivalently, we have
\begin{eqnarray*}
Xu_i=x_iu_i,
\end{eqnarray*}
i.e., $u_i$ is the eigenvector of $M$ with eigenvalue $x_i$. Write
$f=(f_1, \ldots, f_N)$, where $f_i(X)=x_i$, $i=1, \ldots, N$. By
It\^o's formula, we have
\begin{eqnarray}
dx_i(t)=\nabla_{dX_t}f_i(X_t)+{1\over 2}\nabla^2 f_i(X_t)(dX_t,
dX_t).\label{D1}
\end{eqnarray}
By the first order Hadamard variational formula, see \cite{Tao}, we
have
\begin{eqnarray*}
\nabla_{dX_t}f_i(X)&=&u_{i}^*dX_tu_i\\
&=&{1\over \sqrt{N}}u_i^*dB_tu_i-{1\over 2}u_i^*V'(X_t)u_idt.
\end{eqnarray*}
Note that
$$u_i^*V'(X_t)u_i=V'(x_i(t)).$$
By the rotational invariance of the Brownian motion, $U^*B_tU$ is a
Brownian motion on $\mathbb{C}^N$. Denote
$$W_t^i=u_i^*B_tu_i, \ \ \ i=1, \ldots, N.$$
Then
\begin{eqnarray}
\nabla_{dX_t}f_i(X)={1\over \sqrt{N}}dW_t^{i}-{1\over
2}V'(x_i(t))dt.\label{D2}
\end{eqnarray}
On the other hand, by the second order Hadamard variational formula,
see \cite{Tao}, it holds
\begin{eqnarray*}
\nabla^2 f_i(X_t)(dX_t, dX_t)=2\sum\limits_{i\neq j} {|u_{j}^*dX_t
u_i|^2\over x_i(t)-x_j(t)}.
\end{eqnarray*}
By It\^o's calculus, we have
\begin{eqnarray*}
|u_{j}^*dX_t u_i|^2&=&|{1\over \sqrt{N}}u_j^*dB_tu_i-{1\over 2}u_j^*V'(X_t)u_i|^2\\
&=&|{1\over \sqrt{N}}u_j^*dB_tu_i-{1\over 2}V'(x_i(t))u_j^*u_idt|^2\\
&=&|{1\over \sqrt{N}}u_j^*dB_tu_i|^2\\
&=&{1\over N}dt,
\end{eqnarray*}
where we have used the fact that $U^*B_tU$ is a Brownian motion on
$\mathbb{C}^N$. Hence
\begin{eqnarray}
\nabla^2 f_i(X_t)(dX_t, dX_t)={2\over N}\sum\limits_{i\neq j}
{1\over x_i(t)-x_j(t)}dt.\label{D3}
\end{eqnarray}
From $(\ref{D1})$, $(\ref{D2})$ and  $(\ref{D3})$, we derive that
$(x_1(t), \ldots, x_N(t))$ satisfies the following SDE
\begin{eqnarray}
dx_i(t)={1\over \sqrt{N}}dW_t^i-{1\over 2}V'(x_i(t))dt+{1\over
N}\sum\limits_{i\neq j} {1\over x_i(t)-x_j(t)}dt.\label{D4}
\end{eqnarray}
The proof of Theorem \ref{LLX1} is completed. \hfill $\square$

\subsection{From GDBM to matrix diffusions}

The following result provides a random matrix representation for the
generalized Dyson Brownian motion which is defined by solving SDE
$(\ref{SDE1})$.

\begin{theorem}\label{Th-LX2}
Let~$\beta=1, 2$ and~$\lambda_N(0)\in\Delta_N.$~Then, there exists a
$S_N$ (respectively, $\mathcal{H}_N$)-valued diffusion process
$(X^{N, \beta}(t), t\geq 0)$, such that its eigenvalues process
$(\lambda_N(t))_{t\geq0}$  is a solution of the SDE $(\ref{SDE1})$
for the generalized Dyson Brownian motion.
\end{theorem}
{\it Proof}. We only prove Theorem \ref{Th-LX2} for $\beta=1$.
Without loss of generality, we assume that~$X^N(0)$~is the diagonal
matrix D=diag $(\lambda^1_N,\ldots,\lambda^N_N).$ Let~$M>0$~be
fixed. We consider the strong solution $\lambda_N(t)$ of SDE
$(\ref{SDE1})$ untill the stopping time $T_M.$ We let $w_{ij}, 1\leq
i<j\leq N$ be independent Brownian motions. Hereafter, all solutions
will be equipped with the natural filtration $\mathcal
{F}_t=\sigma(w_{ij}(s),W_i(s),s\leq t\wedge T_M)$, where $W_i$ the
Brownian motions of SDE $(\ref{SDE1})$, independent of $w_{ij},
1\leq i<j\leq N.$  For $i<j$, define $R^N_{ij}(t)$ by solving SDE
\begin{eqnarray}
dR^N_{ij}(t)=\sqrt{2\over
N}\frac{1}{\lambda_N^i(t)-\lambda_N^j(t)}dw_{ij}(t), \ \ \
R_{ij}^N(0)=0.\label{ee2.1}
\end{eqnarray}
For $i>j$, set $R_{ij}^N(t)=-R_{ji}^N(t)$. Let $O^N$  be the strong
solution of SDE\footnote{ Here, for two semi-martingales ~$A$ and
$B$~ with values in~$\mathcal {S}_N,$ $\langle
A,B\rangle_t=(\sum\limits^N_{k=1}\langle
A_{ik},B_{kj}\rangle_t)_{1\leq i,j\leq N}$ is the martingale bracket
of $A$ and $B$, and $\langle A\rangle_t$ is the finite variation
part of $A$ at time $t$.}
\begin{equation}\label{E2.3}
\left\{ \begin{aligned}
dO^N(t)&=O^N(t)dR^N(t)-\frac{1}{2}O^N(t)d\langle(R^N)^T,R^N\rangle_t,&\\
O^N(0)&={\rm I}&,
\end{aligned}
\right.
\end{equation} Since SDE $(\ref{E2.3})$ has uniformly
Lipschitz coefficients, we obtain the existence and uniqueness of
strong solutions of $(\ref{E2.3})$ with respect to the filtration
$\mathcal {F}_t$ up to the stopping time $T_M$. By Lemma $4.3.4$ in
\cite{AGZ}, we have $O^N(t)^TO^N(t)={\rm I}$ at all times.

Let $Y^N(t)=O^N(t)^TD(\lambda_N(t))O^N(t)$ and define $W^N(t)$ by
$$dW^N(t)=(O^N(t))^TdY^N(t)O^N(t)$$
with the initial values $W_N(0)=Y_N(0)=X_N(0)$.

By the same argument as used in \cite{AGZ, Gui}, we can prove that
\begin{eqnarray}
dW_N^{ii}(t)=\sqrt{\frac{2}{N}}dW^i_t-\frac{1}{2}
V'(\lambda^i_N(t))dt, \ \ \ \forall i=1, \ldots, N, \label{ee2.6}
\end{eqnarray}and for $i\neq j$,
\begin{eqnarray}
dW_N^{ij}(t)
=\sqrt{2\over N}dw_{ij}(t).\label{ee2.7}
\end{eqnarray}
Since $(O^N(t),t\geq 0)$ is adapted,
$dY^N(t)=O^N(t)dW^N(t)(O^N(t))^T$ is a continuous matrix-valued
semi-martingale. Set
\begin{equation*} \label{eq:1}
 \widetilde{W}_t=\left\{ \begin{aligned}
         &w_{ij}(t), \ \ \ {\rm if}\ i\neq j,\\
                  &W^{i}(t),\ \ \ {\rm if}\ i=j.
                          \end{aligned} \right.
                          \end{equation*}
Then we can write $dY_N(t)=O^N(t)dW^N(t)(O^N(t))^T$ as follows

\begin{eqnarray*}
dY_N(t)&=&O^N(t)\left(\sqrt{2\over N}d\widetilde{W}_t-{1\over 2}D(V'(\lambda_N(t))dt\right)(O^N(t))^T\\
&=&\sqrt{2\over N} O^N(t)d\widetilde{W}_t(O^N(t))^T-{1\over
2}O^N(t)D(V'(\lambda_N(t))(O^N(t))^T dt.
\end{eqnarray*}
Set
\begin{eqnarray*}
B_t=\int_0^t O^N(s)d\widetilde{W}_s (O^N(s))^T.
\end{eqnarray*}
By L\'evy's characterization of Brownian motion,  $B_t$ is a
$S_N$-valued Brownian motion with respect to $\mathcal
{F}_t=\sigma(w_{ij}(s),W_i(s),s\leq t)$. Moreover, by the fact $V$
is a real analytical function and
$Y_N(t)=O^N(t)D(\lambda_N(t))(O^N(t))^T$, we have
\begin{eqnarray*}
\nabla {\rm Tr}V(Y_t)=O^N(t)D(V'(\lambda_N(t))(O^N(t))^T.
\end{eqnarray*}
Hence
\begin{eqnarray}
dY_N(t)=\sqrt{2\over N} dB_t-{1\over 2}\nabla {\rm
Tr}V(Y_t)dt.\label{SDE6}
\end{eqnarray}
Thus, $Y_N(t)$ is the $S_N$-valued diffusion process with generator $L={1\over 2}\Delta_{S_N}-N\nabla{\rm Tr}V\cdot\nabla$. Note that the
SDE $(\ref{SDE6})$  is as the same type as SDE $(\ref{DBM1})$ that $X_N(t)$ satisfies, and $Y_N(0)=X_N(0)$.
 By the uniqueness of
weak solution to SDE, $Y_N(t)$ has the same law as $X_N(t)$, and $\lambda_N(t)$ is the eigenvalue process of $Y_N(t)$. The proof of Theorem \ref{Th-LX2} is completed.\hfill $\square$

\begin{remark} Theorem \ref{Th-LX2} can be considered as a reverse of Theorem \ref{LLX1}.
 The way from $X_N(t)$ to
$\lambda_N(t)$ in Theorem \ref{LLX1} is more direct and does not
need to deal with the question of the existence and uniqueness of
strong solution of the SDE $(\ref{SDE1})$ with singular drift (i.e.,
Theorem \ref{Th1}). On the other hand, the reverse way from
$\lambda_N(t)$ to $Y_N(t)$ in Theorem \ref{Th-LX2} shows that the
generalized Dyson Brownian motion $\lambda_N(t)$ must be obtained by
the way from $X_N(t)$ to $\lambda_N(t)$ as we did in Theorem
\ref{LLX1}. Moreover, Theorem \ref{Th-LX2} also shows that the SDE
$(\ref{SDE1})$ for generalized Dyson Brownian motion has the
uniqueness in distribution. We refer the reader to \cite{No, Gui,
AGZ, Blo} and references therein for further work on matrix-valued
diffusion processes.
\end{remark}

\medskip

\subsection{It\^{o}'s calculus }

 Let $(W^1,\ldots,W^N)$ be
independent Brownian motions and
$(\lambda^1_N(0),\ldots,\lambda^N_N(0))$ be real numbers, let
$\beta\geq 1$ and let $\lambda_N(t)_{t\geq0}$ be the unique strong
solution to SDE $(\ref{SDE1})$. Then by It\^{o}'s calculus, we know
that for all~$f\in {C}^2([0,T]\times\mathbb{R},\mathbb{R}),$~
\begin{eqnarray*}
\int_{\mathbb{R}} f(t,x)L_N(t,dx)&=&\int_{\mathbb{R}}
f(0,x)L_N(0,dx)+\int^t_0\int_{\mathbb{R}}\partial_sf(s,x)L_N(s,dx)ds\\
& &+\frac{1}{2}\int^t_0\int\int_{\mathbb{R}^2}\frac{\partial_xf(s,x)-\partial_yf(s,y)}{x-y}L_N(s,dx)L_N(s,dy)ds\\
& &+\left(\frac{2}{\beta}-1\right)\frac{1}{2N}\int^t_0\int_{\mathbb{R}}\partial_x^2f(s,x)L_N(s,dx)ds\\
& &-\frac{1}{2}\int^t_0\int_{\mathbb{R}}
V'(x)\partial_xf(s,x)L_N(s,dx)ds+M_N^f(t),
\end{eqnarray*}
where $M^f_N$ is the martingale given by
\begin{eqnarray*}
M^f_N(t)=\frac{1}{N}\sqrt{\frac{2}{\beta
N}}\sum\limits^N_{i=1}\int^t_0\partial_xf(s,\lambda^i_N(s))dW^i_s, \ \ \forall\ t\leq T.
\end{eqnarray*}
Note that
\begin{eqnarray*}
\langle M^f_N\rangle_t=\frac{2}{\beta
N^2}\int^t_0\int_{\mathbb{R}}(\partial_xf(s,x))^2L_N(s,dx)ds\leq\frac{2\|\partial_xf\|^2_\infty
t}{\beta N^2}.
\end{eqnarray*}

\section{McKean-Vlasov limit as $N\rightarrow \infty$}

\subsection{Proof of Theorem \ref{Th4} (i)}

We first prove the tightness of $\{L_N(t),t\in [0,T]\}$ in $\mathcal
 {C}([0,T],\mathscr{P}(\mathbb{R}))$,  then we show that the limit of any convergent subsequence of $\{L_N(t),t\in [0,T]\}$ satisfies the
 McKean-Vlasov equation $(\ref{DBM7})$. Here, for $T\in\mathbb{R}^{+},$ we denote by $\mathcal
{C}([0,T],\mathscr{P}(\mathbb{R}))$ the space of continuous
processes from $[0,T]$ into $\mathscr{P}(\mathbb{R})$ equipped with
the weak convergence topology.

\medskip

\textbf{3.1.1 The tightness of the empirical measures}\\

We follow the argument used in \cite{RS93} to  prove the tightness
of $\{L_N(t),t\in[0,T]\}$. Let us pick functions $f_j\in
C_b^\infty(\mathbb{R}, \mathbb{C}), j=1,2,\ldots,$ which is dense in
$C_b(\mathbb{R})$. Thus
$$\langle\mu, f_j\rangle=\langle\mu', f_j\rangle, \ \ \forall
j\Rightarrow  \mu=\mu'.$$ We also pick a $C^\infty$ function
$f_0:\mathbb{R}\rightarrow[1, \infty)$ with the properties
$$f_0(x)=f_0(-x),\ \ \ f_0(x)\rightarrow\infty \ \ \ {\rm as} \ \ x\rightarrow\infty, \ x\in \mathbb{R}^{+}.$$ Taking test functions in the
Schwartz class of smooth functions whose derivatives (up to second order) are rapidly decreasing, we
may assume that
$$f_j, \ f_j'', \ V'f'_j\ \ \ {\rm are\ \ uniformly \ bounded\ \ for\ \ all}\ \
j\geq1.$$ By Ethier and Kurtz \cite{EK} (p.107), to prove the
tightness of $\{L_N(t),\ t\in [0, T],\ N\geq 1\}$,  it is sufficient
to prove that for each $j$ the sequence of continuous real-valued
functions
$$\{\langle L_N(t), f_j\rangle,\ t\in [0, T],\ N\geq1\}$$ is tight.
To this end, note that, by Theorem \ref{Th1}, there is non-collision
and non-explosion for the particles $\lambda_N^i(t)$ for all $t\in
[0, \infty)$. By It\^{o}'s formula in Section 2.2, we have
\begin{eqnarray*}
d\langle L_N(t), f\rangle&=&\frac{1}{N}\sqrt{\frac{2}{\beta
N}}\sum\limits_{i=1}^Nf'(\lambda_N^i(t))dW_t^i+ \left\langle L_N(t),
\left(\frac{2}{\beta}-1\right)\frac{1}{2N}f''-\frac{1}{2}V'f'\right\rangle
dt\\
& &\hskip2cm
+\frac{1}{2}\int\int_{\mathbb{R}^2}\frac{f'(x)-f'(y)}{x-y}L_N(t,dx)L_N(t,dy)dt.
\end{eqnarray*}This yields
\begin{eqnarray}
\langle L_N(t), f_j\rangle&=&\langle L_N(0),
f_j\rangle+\frac{1}{2}\int_0^t\int\int_{\mathbb{R}^2}\frac{f_j'(x)-f_j'(y)}{x-y}L_N(s,dx)L_N(s,dy)ds\nonumber\\
& &-\frac{1}{2}\int_0^t\langle L_N(s), V'f_j'\rangle
ds+\int_0^t\langle L_N(s),
\left(\frac{2}{\beta}-1\right)\frac{1}{2N}f_j''\rangle ds+M_N^{f_j}(t)\nonumber\\
&=&I_1(N)+I_2(N)+I_3(N)+I_4(N)+M_N^{f_j}(t),\label{E3.A}
\end{eqnarray}
where $$M_N^{f_j}(t)=\frac{1}{N}\sqrt{\frac{2}{\beta
N}}\int_0^t\sum\limits_{i=1}^Nf_j'(\lambda_N^i(s))dW_s^i.$$ Note
that, as $L_N(0)$ is weakly convergent, $I_1(N)$ is convergent. By
the assumption that $f_j$ and $f_j''$ are uniformly bounded (hence
$f_j'$ are uniformly bounded) , we can easily show that
$\{M_N^{f_j}(t), t\in [0, T]\}$ and $I_4(N)$ converge  to zero.
Moreover, by the assumption that $V'f_j'$ and $f_j''$ are uniformly
bounded, the Arzela-Ascoli theorem implies that $I_2(N)$ and
$I_3(N)$ are tight in $C([0, T], \mathbb{R})$. Thus the sequence
$\{(L_N(t))_{t\geq0}:\ N\geq1\}$ is tight in $C([0, T],
\mathbb{R})$.  Tightness also follows for $j=0$ if we have
$$\langle L_N(0),f_0\rangle\rightarrow\ \ {\rm finite\ \ limit\ \ as}\ \ N\rightarrow\infty.$$
So let us suppose that the initial distribution $L_N(0)$ have the
property $\langle L_N(0),f_0\rangle\leq K$ for some $K,$ for all
$N.$ For given $\mu_0,$ we could always find $L_N(0)$ and $f_0$ to
satisfy this and the other conditions, and this gives the tightness
for $j=0$ also.

\medskip

\textbf{3.1.2 Identifying of the limit process}\\

Without loss of generality, assuming that $\{L_{N_j}(t), t\in [0,
T]\}$ is a convergent subsequence in $C([0, T],
\mathscr{P}(\mathbb{R}))$. Then,  for all $f\in C_b^2(\mathbb{R})$,
the It\^o's formula $(\ref{E3.A})$ in Section 2.2 and the argument
used in Section 3.1 show that $\langle
\mu(t),f\rangle=\lim\limits_{j\rightarrow \infty}\langle
L_{N_j},f\rangle$ satisfies the following equation
\begin{eqnarray*}
\int_{\mathbb{R}} f(x)\mu_t(dx)&=&\int_{\mathbb{R}}
f(x)\mu_0(dx)+\frac{1}{2}\int_0^t\int\int_{\mathbb{R}^2}\frac{\partial_xf(x)-\partial_yf(y)}{x-y}\mu_s(dx)\mu_s(dy)ds\\
& &\ \ \ \ \ \ \ -\frac{1}{2}\int_0^t\int_{\mathbb{R}}
V'(x)f'(x)\mu_s(dx)ds.
\end{eqnarray*}
That is to say, $\mu(t)$ is a weak solution to the McKean-Vlasov
equation $(\ref{DBM7})$.  Suppose that $\mu_t$
is absolutely continuous with respect to the Lebesgue measure $dx$,
and denote
\begin{eqnarray*}
\rho_t(x)={d\mu_t\over dx}.
\end{eqnarray*}
Integrating by parts in $(\ref{DBM7})$, we can  prove that $\rho_t$
satisfies the nonlinear Fokker-Planck equation $(\ref{NFK1})$. The proof of Theorem \ref{Th4} (i) is
completed.\hfill $\square$

\begin{remark}
To characterize the McKean-Vlasov limit $\mu_t$ of the family of
empirical measures $\{L_N(t), t\in [0, \infty)\}$, we need only to
use the test function $f(x)=(z-x)^{-1}$, where
$z\in\mathbb{C}\backslash\mathbb{R}$, instead of using all test
functions $f\in C_b^2(\mathbb{R})$ in the McKean-Vlasov equation
$(\ref{DBM7})$. Let
\begin{eqnarray*}
G_t(z)=\int_{\mathbb{R}} {\mu_t(dx)\over z-x}
\end{eqnarray*}
be the Cauchy-Stieltjes-Hilbert transform of $\mu_t$. Then $G_t(z)$
satisfies the following equation
\begin{eqnarray*}
{\partial\over \partial t}G_t(z)=-G_t(z) {\partial\over \partial
z}G_t(z)-\frac{1}{2}\int_{\mathbb{R}} {V'(x)\over (z-x)^2}\mu_t(dx).
\end{eqnarray*}
In particular, in the case $V(x)=\theta x^2$, since
\begin{eqnarray*}
-\int_{\mathbb{R}} {x\over (z-x)^2}\mu_t(dx)=z{\partial\over
\partial z}G_t(z)+G_t(z),
\end{eqnarray*}
we obtain
\begin{eqnarray*}
{\partial\over \partial t}G_t(z)=\left(-G_t(z)+\theta z\right)
{\partial\over
\partial z}G_t(z)+\theta G_t(z).
\end{eqnarray*}

\end{remark}

\begin{remark} Following \cite{Gui, AGZ}, we can also prove a stronger version of the tightness of
the empirical measure. More precisely, let $T\in\mathbb{R}^+$, and assuming V satisfy the conditions in
Theorem \ref{Th2}. Assume that
\begin{eqnarray*}
\sup\limits_{N\in \mathbb{N}}\int_{\mathbb{R}}\log(x^2+1)L_N(0, dx)<+\infty. \label{E3.1}
\end{eqnarray*} Then, for all $T\in\mathbb{R}^+$ and $L\in\mathbb{N},$ there exists a compact set $K(L)\subseteq C([0,T],\mathscr{P}(\mathbb{R}))$ such that
\begin{eqnarray*}
\mathbb{P}(L_N(\cdot)\in K(L)^c)\leq e^{-N^2L}.
\end{eqnarray*}
In
particular, the law of $(L_N(s),s\in[0,T])$ is almost surely tight
in $C([0,T],\mathscr{P}(\mathbb{R})).$ This will be used in our forthcoming paper
for a study of the large deviation principle for the generalized
Dyson Brownian motion.  For a proof, see our forthcoming paper.
\end{remark}

\section{McKean-Vlasov equation and optimal transport}

In Section 3.1, we have proved Theorem \ref{Th4} (i), which asserts the existence of weak solution to the McKean-Vlasov equation $(\ref{DBM7})$, equivalently, the existence of weak solution of the nonlinear Fokker-Planck equation $(\ref{NFK1})$. In this section, we use the mass optimal transportation theory to study the uniqueness and the longtime asymptotic behavior of the McKean-Vlasov equation $(\ref{DBM7})$ and the nonlinear Fokker-Planck equation $(\ref{NFK1})$.

Let
\begin{eqnarray*}
W(x)=-2\log |x|, \ \ \ \ x\neq 0.
\end{eqnarray*}
Then the nonlinear Fokker-Planck equation $(\ref{NFK1})$ can be rewritten as follows
\begin{eqnarray}
\partial_t \rho=\nabla\cdot (\rho \nabla (V+W*\rho)).\label{MV-0}
\end{eqnarray}

Before going to study the nonlinear Fokker-Planck equation $(\ref{NFK1})$ (i.e., $(\ref{MV-0})$),  we first recall some
results due to Carrillo, McCann and Villani. In \cite{CMV1},
Carrillo, McCann and Villani studied the McKean-Vlasov evolution
equation of the granular media, which is given by the following
\begin{eqnarray}
\partial_t \rho=\nabla\cdot (\rho \nabla (\log \rho+V+W*\rho)). \label{MV}
\end{eqnarray}
They proved that the McKean-Vlasov equation can be realized as a gradient flow of a free energy functional on the infinite Wasserstein space. More precisely, we have

\begin{theorem} (Carillo-McCann-Villani\cite{CMV1}) Let
\begin{eqnarray}
F(f)=\int_{\mathbb{R}^d} \rho\log \rho dv+\int_{\mathbb{R}^d} \rho V
dv+{1\over 2}\int_{\mathbb{R}^d}\int_{\mathbb{R}^d}
W(x-y)\rho(x)\rho(y)dxdy.\label{F}
\end{eqnarray}
Then the McKean-Vlasov equation $(\ref{MV})$ is the gradient flow of
$F$ with respect to the following infinite dimensional  Riemannian
metric on $\mathscr{P}_2(\mathbb{R}^d)$ (cf. Otto \cite{Ot}):
\begin{eqnarray*}
g_{fdv}(s_1, s_2)=\int_{\mathbb{R}^d}  s_1 s_2 fdv,
\end{eqnarray*}
where $fdv\in \mathcal{P}_2(\mathbb{R}^d)$, $s_1, s_2\in
T_{fdv}\mathscr{P}_2(\mathbb{R}^d)=\{s: \mathbb{R}^d\rightarrow
\mathbb{R}: \int_M s dv=0\}$, and
$$s_i=-\nabla.(f \nabla p_i)$$
for some $p_i\in W^{1,2}(\mathbb{R}^d)$, $i=1, 2$.
\end{theorem}
Moreover, based on Otto's infinite dimensional geometric calculation
on the Wasserstein space, Carrillo, McCann and Villani \cite{CMV1}
proved the following entropy dissipation formula

\begin{theorem}(Carrillo-McCann-Villani\cite{CMV1}) \label{th1} Denote $\xi:=\nabla(\log \rho+V+W*\rho)$. Then
\begin{eqnarray}
{d\over dt} F(\rho_t)&=&-\int_{\mathbb{R}^n} \rho |\xi|^2dv, \label{ED1}\\
{d^2\over dt^2} F(\rho_t)&=&2\int_{\mathbb{R}^n}\rho {\rm Tr}(D\xi)^T(D\xi)dx+2\int_{\mathbb{R}^n} \langle D^2V\cdot \xi, \xi\rangle\rho dx\nonumber\\
& &+\int_{\mathbb{R}^{2n}} \langle D^2W(x-y)\cdot [\xi(x)-\xi(y)],
[\xi(x)-\xi(y)]\rangle d\rho(x)d\rho(y).   \label{ED2}
\end{eqnarray}
\end{theorem}

\medskip

Now let $V: \mathbb{R}\rightarrow [0, \infty)$ be a $C^2$ function
with growth condition $(\ref{grow})$.  In \cite{Bian-Sp01}, Biane
and Speicher proved that the free Fokker-Planck equation
\begin{eqnarray*}
{\partial \rho_t\over \partial t}=-{\partial \over \partial x}(\rho_t({\rm
H}\rho_t-\frac{1}{2}V'))
\end{eqnarray*}
is the gradient flow of  $\Sigma_V$ on the Wasserstein space
$\mathscr{P}_2(\mathbb{R})$. See also \cite{Blo}.

By analogue of the proof of Theorem \ref{th1} in \cite{CMV1}, and observing that for $W(x)=-2\log |x|$ it holds
\begin{eqnarray*}
\xi=&\nabla(V+W*\rho)
=V'-2{\rm H}\rho,
\end{eqnarray*}
we can derive the following dissipation formula for the Voiculescu free
entropy.

\begin{theorem} \label{th2} Let $\xi=V'-2H\rho$. We have
\begin{eqnarray}
\label{entro-diss-3}
{d\over dt}\Sigma_V(\mu_t|\mu_V)&=&-2\int_{\mathbb{R}}\left[V'(x)-2(H\rho)(x)\right]^2\rho(x)dx,\\
{d^2\over dt^2}\Sigma_V(\mu_t|\mu_V)&=&2\int_{\mathbb{R}}  V^{''}(x)|V'(x)-2H\rho(x)|^2\rho(x) dx\nonumber\\
& &+\int\int_{\mathbb{R}^2}
{\left[V'(x)-V'(y)-2(H\rho(x)-H\rho(y))\right]^2\over (x-y)^2}
\rho(x)\rho(y)dxdy.
\end{eqnarray}
\end{theorem}

\medskip

\noindent{\it Proof of Theorem \ref{Th5} (i)}. By Corollary 3.2 in Biane \cite{Bian03}, for any convex $V$, there exists a unique equilibrium measure $\mu$ (indeed $\mu=\mu_V$) with a density $\rho$ satisfying the Euler-Lagrange equation $H\rho(x)={1\over 2}V'(x)$ for all $x\in {\rm supp}(\mu)$.  Thus, $\Sigma_V$ has a unique minimizer $\mu_V$. By the fact that $\Sigma_V$ is lower semi continuous with respect to the weak convergence topology, see e.g. \cite{AGZ, Gui}, we see that it is also lower semi continuous with respect to the Wasserstein topology on $\mathscr{P}(\mathbb{R})$. Moreover,
 as $V$ is $C^2$-convex,  Theorem \ref{th2} implies that
$\Sigma_V$ is a (displacement) convex functional on $\mathscr{P}_2(\mathbb{R})$.

By Proposition $4.1$ in Kloekner \cite{Kloe}, we know that $\mathscr{P}_2(\mathbb{R})$ has vanishing sectional curvature in the sense of Alexandrov. More precisely, for any $\mu_1, \mu_2, \mu_3\in \mathscr{P}_2(\mathbb{R})$ and for any Wasserstein geodesic $\gamma: [0, 1]\rightarrow \mathscr{P}_2(\mathbb{R})$ such that $\gamma(0)=\mu_1$ and $\gamma(1)=\mu_2$, for all $t\in [0, 1]$, it holds that
 \begin{eqnarray*}
 W_2^2(\mu_3, \gamma(t))=tW_2^2(\mu_3, \mu_1)+(1-t)W_2^2(\mu_3, \mu_2)-t(1-t)W_2^2(\mu_1, \mu_2).
 \end{eqnarray*}
Therefore,  $\mathscr{P}_2(\mathbb{R})$ is a nonpositively curved (NPC) space in the sense of Alexandrov (even though $\mathscr{P}_2(\mathbb{R}^n)$ is an Alexander space with nonnegative curvature for $n\geq 2$, see e.g. \cite{AGS}). By Mayer \cite{Mey},  we can conclude that $W_2(\mu_t, \mu_V)\rightarrow 0$ holds if we only assume that $V$ is a $C^2$-convex potential.
The proof of Theorem \ref{Th5} (i) is completed.

\hfill $\square$
\medskip

\noindent{\it Proof of Theorem \ref{Th5} (ii)}.  The proof follows the same
argument as used in \cite{Ot, OV, CMV1}. We use the gradient flow of
the Voiculescu entropy $\Sigma_V$ on the Wasserstein space and the
$K$-convexity along the geodesic displacement between two
probability measures.

Since $V''\geq K$, we have
$$
{d^2\over dt^2}\Sigma_V(\mu_t|\mu_V) \geq K.
$$
By Otto's calculus, we know that
\begin{eqnarray*}
{\rm Hess} \Sigma_V(\mu_t|\mu_V)\left(\frac{\partial \mu_t}{\partial t}, \frac{\partial \mu_t}{\partial t}\right) = {d^2\over dt^2}\Sigma_V(\mu_t|\mu_V),
\end{eqnarray*}
which implies that
$$
{\rm Hess} \Sigma_V(\mu) \geq K.
$$
Let  $\mu(0)=\rho(0)dx$ and $\mu(1)=\rho(1)dx$ be two probability
measures with compact support on $\mathbb{R}$, let
$\mu(s)=\rho(s)dx$ be the unique geodesic in the Wasserstein space
$\mathscr{P}_2(\mathbb{R})$ linking $\mu(0)$ and $\mu(1)$. Then
\begin{eqnarray*}
{d^2\over ds^2}\Sigma_V(\mu(s)) ={\rm Hess}
\Sigma_V(\rho(s))\left(\frac{\partial \rho(s)}{\partial s},
\frac{\partial \rho(s)}{\partial s}\right) \geq
K\left\|\frac{\partial \rho(s)}{\partial
s}\right\|^{2}_{\mathscr{P}_2(\mathbb{R})}.
\end{eqnarray*}
Therefore, for some $\sigma^*\in (0, 1)$,
\begin{eqnarray*}
\Sigma_V(\rho(1)) - \Sigma_V(\rho(0)) &=& \Sigma_V'(\rho(0)) + {1\over 2}\Sigma_V''(\rho(\sigma))|_{\sigma=\sigma^*}\\
&\geq& \left.\left\langle \frac{d\rho(\sigma)}{d\sigma}, \nabla
\Sigma_V \right\rangle\right|_{\sigma = 0} +
{K \over 2}\int^{1}_{0}\left\|\frac{\partial \rho(\sigma)}{\partial \sigma}\right\|^{2}_{\mathscr{P}_2(\mathbb{R})}d\sigma \\
&=&\left.\left \langle \frac{d\rho(\sigma)}{d\sigma}, \nabla \Sigma_V \right\rangle\right|_{\sigma = 0} + {K \over 2} W_{2}^{2}(\rho(0), \rho(1)).\\
\end{eqnarray*}
Similarly,
\begin{eqnarray*}
\Sigma_V(\rho(0)) - \Sigma_V(\rho(1)) &\geq& -\left.\left\langle \frac{d\rho(\sigma)}{d\sigma}, \nabla \Sigma_V \right\rangle\right|_
{\sigma = 1} + {K\over 2} W_{2}^{2}(\rho(0), \rho(1)).\\
\end{eqnarray*}
Summing the two inequalities together, we obtain
\begin{eqnarray}
\left.\left\langle \frac{d\rho(\sigma)}{d\sigma}, \nabla \Sigma_V
\right\rangle\right|_{\sigma = 1} - \left.\left\langle
\frac{d\rho(\sigma)}{d\sigma}, \nabla \Sigma_V
\right\rangle\right|_{\sigma = 0} \geq K W_{2}^{2}(\rho(0),
\rho(1)).\label{DDD1}
\end{eqnarray}

Let $\rho_{t}(s, x)dx:[0, 1] \rightarrow \mathscr{P}_2(\mathbb{R})$
be the unique geodesic between $\mu_t$ and $\mu_V$. By Otto
\cite{Ot}, we have the following  derivative formula of the
Wasserstein distance
\begin{eqnarray}
{d \over dt}W^{2}_{2}(\mu_t, \mu_V) &=& -2\int_{\mathbb{R}}\left.\left\langle\frac{d\rho_t(s)}{ds}(x), \xi_t\right\rangle\right|_{s=0}d\mu_t + 2\int_{\mathbb{R}}\left.\left\langle \frac{d\rho_t(s)}{ds}(x), \xi_t\right\rangle\right|_{s=1} d\mu_V\nonumber\\
&=& 2\left(-\left.\left\langle \frac{d\rho_t(s)}{ds}(x), \nabla \Sigma_V\right\rangle\right|_{s=1} + \left.\left\langle \frac{d\rho_t(s)}{ds}(x),
\nabla \Sigma_V\right\rangle\right|_{s=0}\right)\nonumber\\
&\leq& -2KW_{2}^{2}(\mu_t, \mu_V),\label{WWWW}
\end{eqnarray}
where in the last step we have used $(\ref{DDD1})$. The Gronwall
inequality yields
\begin{eqnarray*}
W_{2}^{2}(\mu_t, \mu_V) \leq e^{-2Kt}W_{2}^{2}(\mu_0, \mu_V).
\end{eqnarray*}

Recall that $$\nabla \Sigma_V(\mu_V)=0.$$  By the fact that $\mu_t$
is the gradient flow of $\Sigma_V$ on $\mathscr{P}_2(\mathbb{R})$
and using the uniform $K$-convexity of $\Sigma_V$, we have
\begin{eqnarray*}
\frac{d}{dt}\|\nabla \Sigma_V(\mu_t)\|^2_{\mathscr{P}_2(\mathbb{R})} &=& 2\left\langle \nabla\|\nabla \Sigma_V(\mu_t)\|^2_{\mathscr{P}_2(\mathbb{R})}, \frac{d\mu_t}{dt}\right\rangle\\
&=& -2Hess \Sigma_V(\mu_t)\left(\frac{d\mu_t}{dt},  \frac{d\mu_t}{dt}\right)\\
&\leq& -2K\left\|\frac{d\mu_t}{dt}\right\|^2_{\mathscr{P}_2(\mathbb{R})}\\
&=& -2K\|\nabla \Sigma_V(\mu_t)\|^2_{\mathscr{P}_2(\mathbb{R})}.
\end{eqnarray*}
Since $\Sigma_V(\mu_V)=0$, we derive that
\begin{eqnarray*}
\frac{d}{dt}\Sigma_V(\mu_t|\mu_V) &=& \left\langle \nabla \Sigma_V(\mu_t), \frac{d\mu_t}{dt}\right\rangle\\
&=& -\|\nabla \Sigma_V(\mu_t)\|^2_{\mathscr{P}_2(\mathbb{R})}\\
&=& \int^{\infty}_{t}\frac{d}{ds}\|\nabla \Sigma_V(\mu_s)\|^2_{\mathscr{P}_2(\mathbb{R})}ds\\
&\leq& -2K\int^{\infty}_{t}\|\nabla \Sigma_V(\mu_s)\|^2_{\mathscr{P}_2(\mathbb{R})}ds\\
&=& 2K\int^{\infty}_{t}\frac{d}{ds}\Sigma_V(\mu_s)ds\\
&=& -2K\Sigma_V(\mu_t|\mu_V),
\end{eqnarray*}
which implies
\begin{eqnarray*}
\Sigma_V(\mu_t|\mu_V)\leq e^{-2Kt}\Sigma_V(\mu_0|\mu_V).
\end{eqnarray*}
The proof of Theorem \ref{Th5} (ii) is completed. \hfill $\square$

\medskip

To prove Theorem \ref{Th5} (iii), we need the following
free logarithmic Sobolev inequality and free Talagrand
transportation cost inequality due to Ledoux and Popescu
\cite{Le-Po09}.

\begin{theorem}\label{LP09}  (Ledoux-Popescu \cite{Le-Po09}) Suppose that $V$ is a $C^2$, convex and there exists a constant $r>0$ such that
$$V''(x)\geq K>0, \ \ \ \ |x|\geq r.$$
Then there exists a constant $c=C(K, r)>0$ such that the free
Log-Sobolev inequality holds: for all probability measure $\mu$ with
$I_V(\mu)<\infty$,
\begin{eqnarray*}
\Sigma_V(\mu|\mu_V)\leq {2\over c}{\rm I}_V(\mu).
\end{eqnarray*}
Moreover, the free Talagrand transportation inequality holds: there
exists a constant $C=C(K, r, V)>0$ such that
\begin{eqnarray*}
CW_2^2(\mu, \mu_V)\leq \Sigma_V(\mu|\mu_V).
\end{eqnarray*}
\end{theorem}

\medskip

\noindent{\it Proof of Theorem \ref{Th5} (iii)}. By Biane and Speicher
\cite{Bian-Sp01},  we have the following entropy dissipation formula
\begin{eqnarray*}
{\partial \over \partial t}\Sigma_V(\mu_t|\mu_V)=-{1\over 2 }{\rm I}_V(\mu_t).\label{entro-diss-2}
\end{eqnarray*}
By Theorem \ref{LP09}, there exists a constant $C_1>0$ such that the {free LSI}
holds
\begin{eqnarray*}
\Sigma_V(\mu|\mu_V)\leq {2\over C_1}{\rm I}_V(\mu),
\end{eqnarray*}
which yields
\begin{eqnarray*}
{d\over dt} \Sigma_V(\mu_t|\mu_V)\leq -{C_1\over 4} \Sigma_V(\mu_t).
\end{eqnarray*}
By the Gronwall inequality, we have
\begin{eqnarray*}
\Sigma_V(\mu_t|\mu_V)\leq e^{-C_1t/4}\Sigma_V(\mu_0|\mu_V).
\end{eqnarray*}
By Theorem \ref{LP09} again, there exists a constant $C_2>0$ such that the  free
transportation cost inequality holds
\begin{eqnarray*}
W_2^2(\mu, \mu_V)\leq {1\over C_2}\Sigma_V(\mu|\mu_V).
\end{eqnarray*}
Therefore
\begin{eqnarray*}
W_2^2(\mu_t, \mu_V)\leq {e^{-C_1t/4}\over C_2}\Sigma_V(\mu_0|\mu_V).
\end{eqnarray*}
This finishes the proof of Theorem \ref{Th4} (iii). \hfill $\square$

\section{Proof of Theorem \ref{Th4} (ii) and (iii)}

{\it Proof of Theorem \ref{Th4} (ii)}. In the proof of Theorem \ref{Th5},
we have proved the following inequalities
\begin{eqnarray*}
\Sigma_V(\mu_2(t))-\Sigma_V(\mu_1(t))\geq \left.\left\langle {\rm
grad}_W \Sigma_V(\rho_t(s)), {\partial\over \partial
s}\rho_t(s)\right\rangle\right|_{s=0} +{K\over 2}W_2^2(\mu_1(t),
\mu_2(t)),
\end{eqnarray*}
and
\begin{eqnarray*}
\Sigma_V(\mu_1(t))-\Sigma_V(\mu_2(t))\geq -\left.\left\langle {\rm
grad}_W \Sigma_V(\rho_t(s)), {\partial\over \partial
s}\rho_t(s)\right\rangle\right|_{s=1} +{K\over 2}W_2^2(\mu_1(t),
\mu_2(t)).
\end{eqnarray*}
Summing them together, we obtain
\begin{eqnarray*}
\left.\left\langle {\rm grad}_W \Sigma_V(\rho_t(s)),
{\partial\over \partial
s}\rho_t(s)\right\rangle\right|_{s=0}-\left.\left\langle {\rm
grad}_W \Sigma_V(\rho_t(s)), {\partial\over \partial
s}\rho_t(s)\right\rangle\right|_{s=1}\leq -{K}W_2^2(\mu_1(t),
\mu_2(t)).
\end{eqnarray*}
By Otto \cite{Ot}, we have the following  derivative formula of the
Wasserstein distance
\begin{eqnarray*}
{d \over dt}W^{2}_{2}(\mu_1(t), \mu_2(t)) &=&
-2\int_{\mathbb{R}}\left.\left\langle\frac{d\rho_t(s)}{ds}(x),
\xi_t\right\rangle\right|_{s=0}d\mu_1(t)
 + 2\int_{\mathbb{R}}\left.\left\langle \frac{d\rho_t(s)}{ds}(x), \xi_t\right\rangle\right|_{s=1} d\mu_2(t)\\
&=& 2\left(-\left.\left\langle \frac{d\rho_t(s)}{ds}(x), \nabla
\Sigma_V(\mu_1(t))\right\rangle\right|_{s=1}
 + \left.\left\langle \frac{d\rho_t(s)}{ds}(x), \nabla \Sigma_V(\mu_2(t))\right\rangle\right|_{s=0}\right)\\
&\leq& -2KW_{2}^{2}(\mu_1(t), \mu_2(t)),
\end{eqnarray*}
which implies
\begin{eqnarray*}
W_{2}(\mu_1(t), \mu_2(t)) \leq e^{-Kt}W_{2}(\mu_1(0),
\mu_2(0)).
\end{eqnarray*}
As a consequence, the McKean-Vlasov equation $(\ref{DBM7})$ has a unique weak solution. This finishes the proof of Theorem \ref{Th4} (ii). \hfill $\square$

\medskip

{\it Proof of Theorem \ref{Th4} (iii)}.
By Theorem \ref{Th4} (i), the family $\{L_N(t), t\in [0, T]\}$ is tight with respect to the weak convergence topology on $\mathscr{P}(\mathbb{R})$, and the limit of any convergent subsequence of  $\{L_N(t), t\in [0, T]\}$ is a weak solution of $(\ref{DBM7})$. By the uniqueness of the weak solution to  $(\ref{DBM7})$, we conclude that $L_N(t)$ weakly converges to $\mu(t)$.
The proof of Theorem \ref{Th4} (iii) is completed.  \hfill $\square$

\medskip

By the same argument as used in Otto \cite{Ot} and Otto-Villani
\cite{OV}, we can prove the following HWI inequalities. To save the
length of the paper, we omit the proof.

\begin{theorem} (HWI inequalities)\ Suppose that there exists a constant $K\in \mathbb{R}$ such that
$$V''(x)\geq K, \ \ \ \forall x\in \mathbb{R}.$$
Let $\mu_i\in \mathscr{P}_2(\mathbb{R})$, $i=1, 2$. Then for all
$t>0$, we have
\begin{eqnarray}
\Sigma_V(\mu_1)-\Sigma_V(\mu_2)\leq W_2(\mu_1, \mu_2)\|{\rm
grad}_W\Sigma_V(\mu_1)\|_{\mathscr{P}_2(\mathbb{R})}-{K\over
2}W_2^2(\mu_1, \mu_2).\label{HWI-1}
\end{eqnarray}
In particular, for any solution to the McKean-Vlasov equation
$(\ref{DBM7})$, we have
\begin{eqnarray}
\Sigma_V(\mu(t))\leq W_2(\mu(t), \mu_V)\|{\rm
grad}_W\Sigma_V(\mu(t))\|_{\mathscr{P}_2(\mathbb{R})}-{K\over
2}W_2^2(\mu(t), \mu(t)).\label{HWI-2}
\end{eqnarray}
where
\begin{eqnarray*}
\|{\rm
grad}_W\Sigma_V(\rho)\|^2_{\mathscr{P}_2(\mathbb{R})}=\int_{\mathbb{R}}\rho
|V'(x)-2H\rho(x)|^2dx.
\end{eqnarray*}
\end{theorem}

\section{Double-well potentials and some conjectures}

Theorem \ref{Th4} and and Theorem \ref{Th5} apply to {$V(x)=a|x|^{p}$} with $a>0$ and $p\geq 2$. When
$a={1\over 2}$, $p=2$ and $\beta=1, 2, 4$, it corresponds to the
GUE, GOE and GSE. Moreover, Theorem \ref{Th4}  and Theorem \ref{Th5} also
apply to the Kontsevich-Penner model on the Hermitian random
matrices ensemble with external potential (cf. \cite{CM})
$$
V(x)={ax^4\over 12}-{bx^2\over 2}-c\log|x|.
$$
Note that, for all $x\neq 0$,
\begin{eqnarray*}
V''(x)&=&ax^2+{c\over x^2}-b\\
 &\geq& {2\sqrt{ac}-b>0}
\end{eqnarray*}
provided that $a>0, c>0$ and
$4ac>b^2$.

Let us consider the double well potential
\begin{eqnarray*}
V(x)={1\over 4}x^4+{c\over 2}x^2, \ \ \ \ x\in \mathbb{R},\label{DW}
\end{eqnarray*}
where $c\in \mathbb{R}$ is a constant. By
\cite{Joh98, BI}, it has been known that the density function of the
equilibrium measure $\mu_V$ is explicitly given as follows:  \\

$(i)$ When $c<-2$,
\begin{eqnarray*}
\rho(x)&=&{1\over 2\pi}|x|\sqrt{(x^2-a^2)(b^2-x^2)}\ \ \ \ \ \ a<|x|<b,\\
&=&0\ \ \ \ \ \ {\rm otherwise},
\end{eqnarray*}
where
$a^2=-2-c$ and $b^2=2-c$. \\

$(ii)$ When $c=-2$, $\rho(x)={1\over 2\pi}x^2\sqrt{4-x^2}$ for $x\in [-2, 2]$ and $\rho(x)=0$ otherwise. \\

$(iii)$ When $c>-2$,
\begin{eqnarray*}
\rho(x)&=&{1\over \pi}(b_2x^2+b_0)\sqrt{a^2-x^2} \ \ \ \ \ \ |x|<a,\\
&=&0\ \ \ \ \ {\rm otherwise},
\end{eqnarray*}
where $a^2={\sqrt{4c^2+48}-2c\over 3}$, $b_0={c+\sqrt{{c^2\over
4}+3}\over 3}$, and $b_2={1\over 2}$.

When $c\in [0, \infty)$, $V$ is $C^2$ convex and $V''(x)\geq 3$ for $|x|\geq 1$. In this case, Theorem \ref{Th4} (iii) implies that $W_2(\mu_t, \mu_V)\rightarrow 0$ with exponential convergence rate.
When $c<-2$, $\mu_V$ has two supports $[-b, -a]$ and $[a, b]$ which are disjoint. By Section $7.1$ in Biane-Speicher \cite{Bian-Sp01}, it is known that $\mu_t$ does not converge to $\mu_V$.  See also Biane \cite{Bian03} and Remark \ref{rem2}. This also indicates that one cannot simultaneously prove a free version of the Holley-Stroock logarithmic Sobolev inequality and a free version of the Talagrand $T_2$-transportation cost
inequality under bounded perturbations of the distribution of eigenvalues $p_N(dx)=Z_N^{-1}\prod_{i<j}|x_i-x_j|^2\prod_{i=1}^Ne^{-NV(x_i)}dx$. Otherwise, by analogue of the proof of Theorem \ref{Th5} (iii), we may prove that $\mu_t$ converges to $\mu_V$ with respect the $W_2$-Wasserstein distance and hence in the weak convergence topology on $\mathscr{P}(\mathbb{R})$. See also \cite{Le-Po09, MMS} for a discussion on non-convex potentials.

In the case $c\in [-2, 0)$, as the global minimizer $\mu_V$ of
$\Sigma_V$ has a unique support, and all stationary point of $\mu_V$
must satisfy the Euler-Lagrange equation $H\mu={1\over 2}V'$, one
can see that the Voiculescu free entropy $\Sigma_V$ has a unique
minimizer $\mu_V$. As $\mu_t$ is the gradient flow of $\Sigma_V$ on
$\mathscr{P}_2(\mathbb{R})$, and since ${d\over
dt}\Sigma_V(\mu_t)=-2\int_{\mathbb{R}}\left[V'(x)-2(H\rho)(x)\right]^2\rho(x)dx$,
we see that $\Sigma_V(\mu_t)$ is strictly decreasing in time $t$
unless $\mu_t$ achieves the (unique) minimizer $\mu_V$. This yields
that $\Sigma_V(\mu_t)$ converges to some value. The question whether
$W_2(\mu_t, \mu_V)\rightarrow 0$ or even $\mu_t$ weakly converges to
$\mu_V$ as $t\rightarrow \infty$ for the above double-well potential
$V$ remains open. We would like to raise the following conjectures.

\begin{conjecture} \label{conj} Consider the double-well potential $V(x)={1\over 4}x^4+{c\over 2}x^2$ with $c\in [-2, 0)$. Then $\mu_t$ converges to $\mu_V$ with respect the $W_2$-Wasserstein distance and hence in the weak convergence topology on $\mathscr{P}(\mathbb{R})$.
  \end{conjecture}

In general, we may raise the following conjectures.

\begin{conjecture}
   Suppose that the potential $V$ is a $C^2$ potential function with $V''(x)\geq K_1$ for all $|x|\geq r$ and $V''(x)\geq -K_2$ for all $|x|\leq r$,
   where $K_1, K_2, r>0$ are some constants. Suppose further that  $\Sigma_V$ has a unique minimizer which has a unique compact support. Then $\mu(t)$ converges  to $\mu_V$ with respect the $W_2$-Wasserstein distance and in the weak convergence topology on $\mathscr{P}(\mathbb{R})$.
\end{conjecture}

\begin{conjecture}
   Suppose that the potential $V$ is a $C^2$ potential function with $V''(x)\geq -K$ for all $x\in \mathbb{R}$, where $K>0$ is a constant.
   Suppose further that  $\Sigma_V$ has a unique minimizer which has a unique compact support. Then $\mu(t)$ converges  to $\mu_V$ with respect the $W_2$-Wasserstein distance and in the weak convergence topology on $\mathscr{P}(\mathbb{R})$.
\end{conjecture}

Finally, let us mention the following conjecture due to Biane and
Speicher \cite{Bian-Sp01}.

\begin{conjecture} Consider the double-well potential given by $V(x)={1\over 2}x^2+{g\over 4}X^4$, where $g$ is a negative constant but very close to zero. Then $\mu_t$ weakly converges to $\mu_V$.
  \end{conjecture}

\begin{flushleft}
\medskip\noindent
Songzi Li, School of Mathematical Science, Fudan University,
220, Handan Road, Shanghai, 200432, China

\medskip

Xiang-Dong Li, Academy of Mathematics and Systems Science, Chinese
Academy of Sciences, 55, Zhongguancun East Road, Beijing, 100190,
China\\
E-mail: xdli@amt.ac.cn
\medskip

Yong-Xiao Xie, Academy of Mathematics and Systems Science,
Chinese Academy of Sciences, 55, Zhongguancun East Road, Beijing,
100190, China
\end{flushleft}


\begin{thebibliography}{99}
\bibitem{AGS} L. Ambrosio, N. Gigli, G. Savr\'e, {\it Gradient Flows
in Metric Spaces and in the Space of Probability Measures}, Birkh\"auser-Verlag, Berlin, 2005.
\bibitem{AGZ} G. Anderson, A. Guionnet, O. Zeitouni, {\it  An Introduction to Random
Matrices}, Cambridge University Press., 2010.
\bibitem{Bian03} P. Biane, Logarithmic Sobolev inequalities, matrix models and free entropy, Acta. Math. Sin. (Engl. Ser.) {\bf 19} (3) (2003), 497-506.
\bibitem{Bian-Sp01} P. Biane, R. Speicher, Free diffusions, free energy and free Fisher information, Ann. Inst. H. Poincar\'e Probab. Stat. {\bf 37} (2001), 581-606.
\bibitem{BCCP1}  D. Benedetto, E. Caglioti, J. A. Carrillo, M. Pulvirenti, A Non-Maxwell Steady distribution for one-dimensional granular media, J. Stat. Phys. 91 (1998), No5/6, 979-990.
\bibitem{BI}P. Bleher, A. Its, Double Scaling Limit in the Random Matrix Model:
The Riemann-Hilbert Approach, Comm. Pure Appl. Math. {\bf 56} (4) (2003), 433-516.
\bibitem{BPS} A. Boutet de Monvel, L. Pastur, M.
Shcherbina, On the statistical mechanics approach in the random
matrix theory: Integrated Density of States, Journal of Statistical
Physics. Vol. 79. Nos. 3/4. 1995.
\bibitem{Blo} G. Blower, {\it Random Matrices: High Dimensional
Phenomena}, London Math. Soc., Lect. Notes Ser. 367, Cambridge Univ.
Press, 2009.
\bibitem{BZ93} E. Brezin, A. Zee, Universality of the corrections between
eigenvalues of large random matrices, Nucl. Phys. B 402, 613-627
(1993).
\bibitem{CM} L. Chekhov, Yu. Makeenko, The multicritical Kontsevich-Penner model, Morden Phys. Letters A. Vol. {\bf 7} No. 14 (1992), 1223-1236.
\bibitem{CMV1}J. Carrillo, R. McCann, C. Villani, Kinetic equilibration rates for granular media and related equations: entropy dissipation and mass transportation estimates, Rev. Mat. Iberoamericana {\bf 19} (2003), 971-1018.
\bibitem{Ch} T. Chan, The Wigner semi-circle law and eigenvalues of
matrix-valued diffusions. Probab.Theory Relat.Fields 93, 249-272
(1992).
\bibitem{Deift} P. Deift,  {\it Orthogonal Polynomials and Random Matrices: A Riemann-Hilbert
Approach}, American Mathematical Society, 2000.
\bibitem{Dy1}  F. J. Dyson, Statistical theory of the energy
levels of complex systems. I, II, and III. J. Math. Phys. 3.
140-156, 157-165, 166-175, 1962.
\bibitem{Dy2} F. J. Dyson,  A Brownian-motion model of the eigenvalues of
a random matrix. J. Mathematical Phys. 3. 1191-1198, 1962.
\bibitem{EK} S.N. Ethier, T.G. Kurtz, {\it Markov Processes:
Characterization and convergence}, New York, Wiley 1986.
\bibitem{FFS92} R. Fernandez, J. Fr\"{o}hlich, A. Sokal, {\it Random Walks,
Critical Phenomena, and Triviality in Quantum Field Theory},
Springer-Verlag, Heidelberg, 1992.
\bibitem{Gui} A. Guionnet,  {\it Large random matrices: Lectures
on macroscopic asymptotics}, Springer, 2008.
\bibitem{IW} N. Ikeda, S. Watanabe, {\it Stochastic
Differential Equations and Diffusion Processes}, North-Holland
Publishing Company, Amsterdam, 1981.
\bibitem{Joh98} K. Johansson, On fluctions of eigenvalues of random Hermitian matrices, Duke Math. J. {\bf 91} (1998), 151-204.
\bibitem{Kon} M. Kontsevich, Intersection theory on the moduli space of curves
and the matrix Airy function, Commun. Math. Phys. 147 (1992), 1-23.
\bibitem{KT} M. Katori, H. Tanemura, Non-equilibrium dysnamics of Dyson's model with an infinite number of particles, preprint, 2009.
\bibitem{Lan} C. Lancelloti, On the fluctuations about the Vlasov limit for $N$-particle systems with mean-field interactions, J. Stat. Phys. {\bf 136} (2009), 643-665.
\bibitem{Kloe} B. Kloeckner, A geometric study of Wasserstein spaces: Euclidean spaces, Annali della Scuola Normale Superiore di Pisa, Classe di Scienze IX, {\bf 2} (2010) 297-323, DOI10.2422/2036-2145.2010.2.03, and www.arXiv:0804.3505.
\bibitem{Le-Po09} M. Ledoux, I. Popescu, Mass transportation proofs of free functional inequalities, and free Poincar\'e inequalities, J. Funct. Anal. {\bf 257} (2009), 1175-1221.
\bibitem{MMS}M. Ma\"{i}da, E. Maurel-Segala, Free transport-entropy
inequalities for non-convex potentials and application to
concentration for random matrices, hal-00687686 version 2, September
10, 2012.
\bibitem{Mal} F. Malrieu, Convergence to the equilibrium for
granular media equations and their Euler schemes, Ann. Appl. Probab.
13 (2003), No. 2, 540-560.
\bibitem{Meh} M. Mehta, {\it Random Matrices}, 3rd Edition, Elsevier (Singapore) Pte
Ltd, 2006.
\bibitem{Mey}U.F. Mayer, Gradient flows on non positively curved metric
spaces and harmonic maps, Commun. in Analysis and Geome. {\bf 6}, (1998), No. 2, 199-253.
\bibitem{No}J. R. Norris, L.C.G. Rogers, D. Williams, Brownian
motions of ellipsoids, Trans. Amer. Math. Soc. 294, 2 (1986),
757-765.
\bibitem{Os}H. Osada,  Interacting Brownian motions in infinite dimensions with
logarithmic interaction potentials, arXiv:0902.3561v2
\bibitem{Ot}  F. Otto, The geometry of dissipative evolution equations: the porous medium equation, Commun.  Parial Differential Equations 26 (2001), No.1/2,  101-174.
\bibitem{OV} F. Otto, C. Villani, Generalization of an inequality by Talagrand and links with the logarithmic Sobolev
inequality, J. Funct. Anal. 173 (2000) 361-400.
\bibitem{PaSh} L. Pastur, M. Shcherbina, {\it Eigenvalue Distribution of Large Random Matrices}, Amer. Math. Soc., Mathematical Surveys and Monographs, Vol 171, 2010.
\bibitem{Ro} C. Prevot, M. Rockner,  {\it A Concise Course on Stochastic Partial Differential
Equations}, Lecture Notes in Math. 1905, Springer-Verlag, Berlin,
2007.
\bibitem{RS93}  L.C.G. Rogers, Z. Shi, Interacting
brownian particles and the Wigner law, Probab. Theory Related Fields
95,4(1993), 555-570.
\bibitem{Sph} H. Spohn, {\it Large Scale Dynamics of Interacting Particles}, Springer-Verlag, Berlin Heidelberg 1991.
\bibitem{Tao} T. Tao, {\it Topics in Random Matrices Theory}, Graduate Studies in Mathematics, Vol. 132, Amer. Math. Soc. 2012.
\bibitem{Vi1} C. Villani, {\it Topics in Mass Transportation}, Grad. Stud. Math., Amer. Math. Soc., Providence, RI, 2003.
\bibitem{Vi2} C. Villani, {\it Optimal Transport, Old and New}, Springer-Verlage, Berlin, 2009.
\bibitem{Voi1} D. Voiculescu, The analogues of entropy and of Fisher's information
measure in free probability theory, I, Commun. Math. Phys. 155
(1993), 71-92.
\bibitem{W} E. P. Wigner, On the distribution of the roots of
certain symmetric matrices. Ann. of Math. (2) {\bf 67} (1958), 325-327.
\end{thebibliography}
\end{document}